\documentclass[journal]{IEEEtran}
\usepackage{amsfonts}
\usepackage{bm}
\usepackage{amsmath,amsfonts,amsthm,amssymb}
\usepackage{setspace}
\usepackage{Tabbing}
\usepackage{fancyhdr}
\usepackage{lastpage}
\usepackage{extramarks}
\usepackage{chngpage}
\usepackage{algorithmic}
\usepackage{soul,color}
\usepackage{epsfig}
\usepackage{graphicx,float,wrapfig,subfigure}
\usepackage{dsfont}
\usepackage{longtable}
\usepackage{wrapfig}
\usepackage{stfloats}
\usepackage{cite}
\usepackage{graphicx}
\usepackage{array}
\usepackage{multirow}{\tiny }
\usepackage{multicol}
\usepackage{colortbl}
\usepackage{tabularx}
\usepackage{mdwmath}
\usepackage{mdwtab}
\usepackage{color}
\usepackage{verbatim}
\usepackage{amsmath}
\usepackage{tikz}
\usetikzlibrary{calc}
\usepackage{flowchart}
\usetikzlibrary{shapes.geometric, arrows}
\usepackage{overpic}
\usepackage{epstopdf}
\usepackage{url}
\usepackage[colorlinks,linkcolor=black,anchorcolor=black,citecolor=black,urlcolor=black]{hyperref} 

\usepackage{makecell}
\usepackage{textcomp}
\usepackage{algorithm}
\usepackage{algorithmic}

\newtheorem{remark}{Remark}

\usepackage{threeparttable}

\usepackage{color}
\makeatletter %
\let\myorg@bibitem\bibitem
\def\bibitem#1#2\par{%
	\@ifundefined{bibitem@#1}{%
		\myorg@bibitem{#1}#2\par
	}{%
		\begingroup
		\color{\csname bibitem@#1\endcsname}%
		\myorg@bibitem{#1}#2\par
		\endgroup
	}%
}

\makeatother %

\allowdisplaybreaks

\usepackage{footnote}
\makesavenoteenv{tabular}
\makesavenoteenv{table}

\begin{document}

\title{
Efficient constraint learning for data-driven active distribution network operation  
	}

\author{
Ge~Chen,~\IEEEmembership{Graduate Student Member,~IEEE,}
Hongcai~Zhang,~\IEEEmembership{Member,~IEEE,}
and~Yonghua~Song,~\IEEEmembership{Fellow,~IEEE}
\vspace{-6mm}
}

\maketitle

\begin{abstract}
Scheduling flexible sources to promote the integration of renewable generation is one fundamental problem for operating active distribution networks (ADNs). However, existing works are usually based on power flow models, which require network parameters (e.g., topology and line impedance) that may be unavailable in practice. To address this issue, we propose an efficient constraint learning method to operate ADNs. This method first trains multilayer perceptrons (MLPs) based on historical data to learn the mappings from decisions to constraint violations and power loss. Then, power flow constraints can be replicated by these MLPs without network parameters. We further prove that MLPs learn constraints by formulating a union of disjoint polytopes to approximate the corresponding feasible region. Thus, the proposed method can be interpreted as a piecewise linearization method, which also explains its desirable ability to replicate complex constraints. Finally, a two-step simplification method is developed to reduce its computational burden. The first step prunes away unnecessary polytopes from the union above. The second step drops redundant linear constraints for each retained polytope. Numerical experiments based on the IEEE 33- and 123-bus test systems validate that the proposed method can achieve desirable optimality and feasibility simultaneously with guaranteed computational efficiency.
\end{abstract}

\begin{IEEEkeywords}
Active distribution networks, deep learning, renewable generation, optimal power flow, flexible sources.
\end{IEEEkeywords}

\section{Introduction} \label{sec_intro}
\IEEEPARstart{M}{otivated} by the carbon neutrality goal, increasing renewable generation (DRG) has been integrated into distribution networks \cite{heptonstall2021systematic}. The high penetration of DRG promotes the transformation of traditional passive distribution networks into active distribution networks (ADNs), which unlocks the potential flexibility from the distribution level \cite{9151405}.

For ADNs, one fundamental problem is how to coordinate their flexible sources, such as energy storage systems \cite{8463580}, electric vehicles \cite{8758383}, or heating, ventilation, and air conditioning (HVAC) systems \cite{9417102}, to reduce the curtailment of DRG. Traditionally, this coordination is guided by power flow models. For example, reference \cite{9372881} employed the second-order cone programming (SOCP) relaxation of the DistFlow model to coordinate heating systems with DRG.
Based on this SOCP relaxation, reference \cite{9112232} designed a distributed algorithm for managing the demand responses of HVAC systems in ADNs. 
Reference \cite{9628048} also used this SOCP relaxation to design the charging/discharging scheduling strategy for electric vehicles to promote DRG integration in ADNs. Reference \cite{9314264} combined a linearized power flow model with distributed optimization schema to schedule flexible sources in ADNs. In general, these power flow-based works can effectively operate ADNs. However, they usually require power network parameters (e.g., network topology and line impedance) that are often unavailable at the distribution level due to unaware topology changes or inaccurate data maintenance \cite{7875102}.
Moreover, some approximations (e.g., SOCP relaxations used in \cite{9372881,9112232,9628048}) are usually introduced to convexity the non-convex power flow equations.
However, these approximations may introduce considerable errors and even infeasible solutions \cite{6815671}. 

Due to the popularization of smart meters, collecting operational data of ADNs (e.g., nodal power injections, bus voltages, and power flows) has become easier and cheaper \cite{8322199}.
Since this data contains the knowledge of network modeling, learning-based methods may be an alternative choice for operating ADNs \cite{9265482}. Generally speaking, these learning-based methods can be divided into three categories.

\subsubsection{Optimize-then-learn methods}
Methods in this category first solve multiple instances of optimization problems to construct training sets, where the input features are operating conditions (e.g., power demands and available DRG) and the output labels are optimal solutions (e.g., schedules of flexible sources), respectively. Then, supervised learning models are trained to approximate the mapping from a specific condition to the optimal decision. 
For example, reference \cite{9205647} trained multilayer perceptrons (MLPs) to predict the optimal generation based on demand conditions. References \cite{chatzos2020high,9335481} combined this method with the Lagrangian dual approach to improve its feasibility. Reference \cite{9314196} replaced MLPs with recurrent neural networks to predict the optimal schedule of flexible sources. Since the above methods replace the solving process of optimization problems with neural networks' inferences, they can achieve excellent computational efficiency.

\subsubsection{Reinforcement learning}
Reinforcement learning (RL) trains agents on how to act to maximize the expected cumulative reward (e.g., the opposite of total operation costs) in the future. Reference \cite{8944292} used RL to realize the voltage control of distribution networks. Reference \cite{9353702} leveraged the multi-agent RL to operate distribution networks, which can reduce communication costs required by network training. Reference \cite{9508140} proposed a graph-based RL method to coordinate flexible sources for the restoration of distributed systems. Since RL can directly infer the best action based on current state information, it shows great potential in realizing the real-time operation of ADNs.

\subsubsection{Constraint learning}
Constraint learning trains neural networks to learn and replicate complex constraints, e.g., non-convex power flow equations. For example, in references \cite{9302963,venzke2020neural}, binary classifiers were trained as surrogate models of power flow models to judge the feasibility of a given decision. After training, the classifiers were equivalently reformulated as mixed-integer linear constraints that off-the-shelf solvers can easily handle. 
Instead of using binary classifiers, our previous work employed a regression multi-layer perceptron (MLP) to predict the maximum constraint violation so that the feasibility of solutions can be improved \cite{9502573}. By replacing the MLP with a quantile regression network, we further extended constraint learning to chance-constrained scheduling problems for distribution networks \cite{chen2022deep}. Constraint learning only requires operational data instead of optimal solutions as training data, so it can bypass the requirement of network parameters. Moreover, desirable optimality can also be achieved since the global optima of the mixed-integer replication can be found by the Branch-and-Bound algorithm.

Although the above papers have confirmed the effectiveness of the learning-based methods, they are still facing some challenges. Specifically, the optimize-then-learn methods, including \cite{9205647,chatzos2020high,9335481,9314196}, need optimal solutions as training labels. Thus, network parameters are still indispensable because they need to build power flow models to generate training labels. For RL-based methods, such as \cite{8944292,9353702,9508140}, their agents need to continuously interact with actual environments (e.g., distribution networks) to update 
their policies. This ``trial and error" scheme can be risky and is usually unacceptable in practice. The constraint learning methods used in \cite{9302963,venzke2020neural,9502573,chen2022deep} do not require network parameters or interactions with real systems. However, since neural networks are usually regarded as ``black boxes", there is a lack of interpretation to explain why constraint learning can replicate complex constraints. Meanwhile, the computational efficiency of constraint learning may be undesirable because of its mixed-integer characteristic. Moreover, reference \cite{yang2021modeling} pointed out that constraint learning may extrapolate significantly from the training data, leading to a solution in an area where neural networks have not learned. In this case, the 
reliability of solutions can not be ensured.

On the basis of our previous works \cite{9502573,chen2022deep}, we propose an interpretable and more efficient constraint learning method to overcome the challenges above. The proposed method trains two MLPs based on historical samples to predict constraint violations and power loss. Then, power flow constraints can be replicated with no need for network parameters. Compared to our previous works and other existing methods, the main contributions of this paper are twofold:

\begin{enumerate}
\item We prove that constraint learning can be interpreted as a piecewise linearization (PWL) method for multivariate functions. 
Specifically, the trained MLP first partitions the domain of decision variables (i.e., schedules of flexible sources) into multiple disjoint regions. In each region, the MLP's output (i.e., constraint violation) is linearly dependent on the decision variables. By requiring the predicted violation to be non-positive, the feasible region of the original problem can be approximated by a subset of the space covered by these disjoint regions. We further prove that this subset can be represented as a union of polytopes. Since the polytope number is positively correlated with the MLP's neuron number, the accuracy of constraint learning can be guaranteed once the MLP has sufficient neurons. 
\item 
Considering that constraint learning may introduce too many polytopes and lead to a huge computational burden, we develop a two-step simplification method, resulting in a more efficient formulation. In the first step, we prune away unnecessary polytopes. In the second step, we drop redundant linear constraints for each retained polytope. Then, the computational burden can be reduced. Moreover, this simplification method does not allow significant extrapolation from historical data. Thus, the reliability of solutions can also be enhanced.
\end{enumerate}

The remaining parts are organized as follows. Section \ref{sec_modeling} describes the problem formulation of ADNs' operation. Section \ref{sec_solution} introduces the proposed PWL-based interpretation and simplification method. Section \ref{sec_case} demonstrates simulation results, and Section \ref{sec_conclusion} concludes this paper.

\section{Problem formulation} \label{sec_modeling}
We focus on a radial ADN with flexible resources. The operating goal is to minimize the total cost while satisfying all critical constraints. For simplicity, this paper uses HVAC systems as an example to represent flexible resources. Note that other flexible resources can also be considered in the same framework by adjusting the expression of nodal power injections. In this paper, we use non-Bold and Bold lowercase letters to represent scalars and vectors, while matrices are denoted by Bold uppercase letters.

\subsection{System modeling}
\subsubsection{HVAC systems}
Due to buildings' inherent ability to store heating/cooling power, the power profile of HVAC systems can be adjusted to promote the integration of DRG while keeping indoor thermal comforts.
By using $i \in \mathcal{I}$ and $t \in \mathcal{T}$ to index HVAC systems and time slots, the thermodynamic model of one building can be expressed as \cite{8731716}:
\begin{align}
\theta_{i,t}^\text{in}=&a_i^\text{in}\theta_{i,t-1}^\text{in}+a_i^\text{out}\theta_{t-1}^\text{out} \notag \\
&+ a_i^\text{h}(q_{i,t-1}^\text{heat} - q_{i, t-1}^\text{cool}), \ \forall i \in \mathcal{I}, \ \forall t \in \mathcal{T}, \label{eqn_thermal}
\end{align}
where $\theta^\text{in}_{i, t}$ and $\theta^\text{out}_t$ are the indoor and outdoor temperatures, respectively; $q_{i,t-1}^\text{heat}$ and $q_{i, t-1}^\text{cool}$ denote the heat load contributed by indoor sources (e.g., humans and electronic devices) and cooling supply of HVAC systems, respectively. Parameters $a_i^\text{in}$, $a_i^\text{out}$, and $a_i^\text{h}$ are calculated by:
\begin{align}
a_i^\text{in}=e^{-\frac{g_i}{C_i}\Delta t}, \  a_i^\text{out}=1-a_i^\text{in}, \ 
a_i^\text{h}=a_i^\text{out}/g_i, \ \forall i \in \mathcal{I},
\end{align}
where $C_i$ and $g_i$ are the building heat capacity and heat transfer coefficient between indoor and outdoor environments, respectively; $\Delta t$ is the length of a time interval. Indoor thermal comforts require that all indoor temperatures should maintain in comfortable regions:
\begin{align}
\bm \theta^\text{min} \leq \bm \theta_{t}^\text{in} \leq \bm \theta^\text{max}, \forall t \in \mathcal{T},
\end{align}
where $\bm \theta^\text{min}$ and $\bm \theta^\text{max}$ are the lower and upper bounds of comfortable regions, respectively. 
The active and reactive powers of HVAC systems are modeled by:
\begin{align}
p_{i,t}^\text{HV} = \frac{q_{i,t}^\text{cool}}{\text{COP}_i}, \ q_{i,t}^\text{HV} = \frac{\sqrt{1 - \phi_i^2}}{\phi_i}  p_{i,t}^\text{HV}, \ \forall i \in \mathcal{I}, \forall t \in \mathcal{T},
\end{align}
where $\text{COP}_i$ is the coefficient of performance of the $i$-th HVAC system; $\phi_i$ is the power factor of the $i$-th HVAC system. The active power of HVAC systems should also be kept in allowable regions due to device limitations:
\begin{align}
\bm p_{t}^\text{HV} \leq \bm p_{t}^\text{HV,max}, \quad \forall t \in \mathcal{T}.
\end{align}

\subsubsection{Distributed renewable generation}
The actual outputs of distributed renewable generators $\bm p_t^\text{DG}$ can be controlled by adjusting their curtailment rate $\bm \lambda_t$, as follows:
\begin{align}
\bm p_t^\text{DG} = \bm G_{t}^\text{DG} * (\bm 1-\bm \lambda_t), \  \bm 0 \leq \bm \lambda_t \leq \bm 1, \ \forall t \in \mathcal{T}, \label{eqn_DG}
\end{align}
where $\bm G_{t}^\text{DG}$ is the available DRG; operator $*$ represents the element-wise multiplication.

\subsubsection{Nodal power injections}
The power injections on each node (except the root node), i.e., $\bm p_t$ and $\bm q_t$, are expressed as:
\begin{align}
\bm p_t = -\bm p_t^\text{HV} - \bm p_{t}^\text{d}  + \bm p_t^\text{DG}, \ 
\bm q_t = -\bm q_t^\text{HV} - \bm q_{t}^\text{d}, \ \forall t \in \mathcal{T}, \label{eqn_q}
\end{align}
where $\bm p_{t}^\text{d}$ and $\bm q_{t}^\text{d}$ are the base active and reactive demands, i.e., the demands except HVAC loads. 

\subsubsection{Power flow constraints}
The power flow of a radial ADN can be described by the DistFlow model \cite{19266}:
\begin{align}
\begin{cases}
\sum_{k \in \mathcal{C}_{j}} P_{jk,t} = p_{j,t} + P_{ij,t} - r_{ij}I_{ij,t}^2, \\
\sum_{k \in \mathcal{C}_{j}} Q_{jk,t} = q_{j,t} + Q_{ij,t} - x_{ij}I_{ij,t}^2, \\
V_{j,t}^2=V_{i,t}^2-2(r_{ij}P_{ij,t}+x_{ij}Q_{ij,t})\\
\quad \quad \quad \quad \quad+ (r_{ij}^2 + x_{ij}^2)I_{ij,t}^2, \\
I_{ij,t}^2 = \frac{P_{ij,t}^2 + Q_{ij,t}^2}{V_{i,t}^2},\\
\forall (i,j) \in \mathcal{B},\forall t \in \mathcal{T},
\end{cases}\label{eqn_distflow}
\end{align}
where $P_{ij,t}$ and $Q_{ij,t}$ are the active and reactive power flows on line $(i, j)$, respectively; $p_{j,t}$ and $q_{j,t}$ denote the active and reactive power injections on bus $j$, respectively; $V_{i,t}$ and $I_{ij,t}$ are the voltage and current magnitudes on bus $i$ and line $(i, j)$, respectively; $r_{ij}$ and $x_{ij}$ denotes the resistance and reactance of line $(i,j)$, respectively. Symbol $\mathcal{B}$ represents the index set of lines. Set $\mathcal{C}_{j}$ contains the child bus indexes of bus $j$. 
To guarantee security, the magnitudes of bus voltages and power flows shall stay in specific regions:
\begin{align}
\underline{\bm V} \leq \bm V_t \leq \overline{\bm V}, \ 
\sqrt{\bm P_t^2 + \bm Q_t^2} = \bm S_t \leq \overline{\bm S}, \ \forall t \in \mathcal{T}, \label{eqn_OPF_constraint}
\end{align}
where $\underline{\bm V}$ and $\overline{\bm V}$ are the lower and upper bounds of voltage magnitudes; $\bm S_t$ and $\overline{\bm S}$ are the actual and maximum allowable values of apparent power flows on lines.

\subsubsection{Total cost}

The energy purchasing cost in one time slot, $EC_t$, can be calculated by:
\begin{align}
	&EC_t =  (\eta^\text{buy} G_t^\text{buy} - \eta^\text{sell} G_t^\text{sell})\Delta t,\  \forall t \in \mathcal{T}, \label{eqn_EC}\\
	&G_t^\text{buy} - G_t^\text{sell} = G_t^\text{root}, \  G_t^\text{buy}\geq0,\  G_t^\text{sell}\geq0,\  \forall t \in \mathcal{T}, \label{eqn_netload}
\end{align} 
where $\eta^\text{buy}$ and $\eta^\text{sell}$ represent the per-unit prices of electricity purchasing and selling, $\eta^\text{buy}\geq \eta^\text{sell}$; $G_t^\text{buy}$ and $G_t^\text{sell}$ are two auxiliary variables. Variable $G_t^\text{root}$ is the net active power at the root node, which can be calculated by the network-level power balance:
\begin{align}
G_t^\text{root} = \bm 1^\intercal \bm p_t + p_t^\text{loss}, \  \forall t \in \mathcal{T}, \label{eqn_balance} 
\end{align}
where $p_t^\text{loss}$ is the total power loss and can be calculated by:
\begin{align}
p_t^\text{loss} = \sum_{(i,j) \in \mathcal{B}} I_{ij,t}^2 r_{ij}, \  \forall t \in \mathcal{T}. \label{eqn_loss}
\end{align}

\subsection{Formulation of the optimization problem}
The operating goal for an ADN is to minimize the energy purchasing cost while meeting all critical constraints, which can be summarized as \textbf{P1}:
\begin{align} 
	&\min_{(\bm p_t^\text{HV}, \bm \lambda_t)_{\forall t \in \mathcal{T}}} \quad \sum_{t \in \mathcal{T}} EC_t, &\text{s.t.: Eqs.} &\text{ (\ref{eqn_thermal})-(\ref{eqn_loss})} \tag{$\textbf{P1}$}.
\end{align}
As mentioned in Section \ref{sec_intro}, formulating \textbf{P1} is challenging because it requires network parameters, which may be unavailable in practice. Moreover, even if all parameters are known, \textbf{P1} is hard to solve due to the non-convex constraint (\ref{eqn_distflow}).

\section{Solution Methodology} \label{sec_solution}
We propose an efficient constraint learning method to overcome the above challenges. In this section, we first introduce the conventional constraint learning in detail. Then, a PWL-based interpretation is proposed to explain why constraint learning can replicate constraints. Finally, a simplification method is developed to improve its computational efficiency.

\subsection{Introduction of constraint learning} \label{sec_CL}
The key idea of constraint learning is to replace complex constraints with trained neural networks. Following reference \cite{9502573}, two MLPs with ReLU as activation functions are trained to replicate power flow constraints. The first MLP predicts constraint violations (termed as ``Vio-MLP"), while the second one forecasts the total power loss (termed as ``Loss-MLP"). 

\subsubsection{Vio-MLP}
The input of the Vio-MLP, i.e., $\bm x_t$, is defined as the collection of the active/inactive power demands and the actually used DRG, as follows:
\begin{align}
&\bm x_t = (\underbrace{\bm p_t^\text{HV}+\bm p_t^\text{d}}_\text{active demand},\  \underbrace{\bm q_t^\text{HV}+\bm q_t^\text{d}}_\text{reactive demand},\  \underbrace{\bm p_t^\text{DG})}_\text{used DRG},\ \forall t \in \mathcal{T}. \label{eqn_x}
\end{align}
The output of the Vio-MLP, i.e., $h_t$, is defined as a measurement of power flow constraint violations, as follows:
\begin{align}
&h_t = \max \{\frac{\underline{\bm V} -  \bm V_t}{\overline{\bm V}-\underline{\bm V}}, \ \frac{\bm V_t - \overline{\bm V}}{\overline{\bm V}-\underline{\bm V}}, \ \frac{\bm S_t - \overline{\bm S}}{\overline{\bm S}}\}, \ \forall t \in \mathcal{T}, \label{eqn_h}
\end{align}
where $h_t$ is designed based on (\ref{eqn_OPF_constraint}). Note that we standardize the violations of (\ref{eqn_OPF_constraint}) to make the violations of voltage and power flow limitations comparable. Based on the operational data of the ADN, the mapping from $\bm x_t$ to $h_t$ can be approximately represented by the forward propagation of a trained MLP with ReLU as the activation function:
\begin{align}
&\bm s_t^{0} = \bm x_t, \ \forall t \in \mathcal{T}, \label{eqn_input}\\
&\bm z_t^{l} = \bm W^{l} \bm s_t^{l-1} + \bm b^{l}, 
\forall l \in \mathcal{L}, \ \forall t \in \mathcal{T}, \label{eqn_z} \\
&\bm s_t^{l} = \max(\bm z_t^{l}, 0), \quad \forall l \in \mathcal{L}, \ \forall t \in \mathcal{T}, \label{eqn_s} \\
&h_t = (\bm w^{|\mathcal{L}|+1})^\intercal \bm s_t^{|\mathcal{L}|} + b^{|\mathcal{L}|+1},\ \forall t \in \mathcal{T},\label{eqn_output}
\end{align}
where $\bm z^{l}$ and $\bm s^{l}$ are the outputs of the linear mapping and activation function in hidden layer $l$; $\mathcal{L}$ denotes the set of hidden layers, $l \in \mathcal{L}$; parameters ($\bm W^l, \bm b^{l})_{\forall l \in \mathcal{L}}$ and $(\bm w^{|\mathcal{L}|+1}, b^{|\mathcal{L}|+1})$ are the weights and bias that are to be learned through training.

By restricting the measurement $h_t$ to be non-positive:
\begin{align}
h_t \leq 0, \ \forall t \in \mathcal{T}, \label{eqn_safe}
\end{align}
the feasibility of the power flow constraints are guaranteed. Hence, the non-convex power flow constraints (\ref{eqn_distflow})-(\ref{eqn_OPF_constraint}) can be replicated by (\ref{eqn_input})-(\ref{eqn_safe}). Considering that the maximum operator in (\ref{eqn_s}) is intractable for off-the-shelf solvers, reference \cite{9502573} further equivalently reformulated constraints (\ref{eqn_z})-(\ref{eqn_s}) into a solvable mixed-integer linear form, as follows:
\begin{align}
	&\left\{
	\begin{aligned}
	&\bm s_t^{l} - \bm r_t^{l}=\bm W^{l} \bm s_t^{l-1} + \bm b^{l}, \\
	&0 \leq \bm s_t^{l} \leq M \cdot \bm \mu_t^{l},\\
	&0 \leq \bm r_t^{l} \leq M\cdot(1-\bm \mu_t^{l}),\\
	&\bm \mu_t^{l} \in \{0,1\}^{N_l},
	 \end{aligned}\right. \ \forall l \in \mathcal{L}, \ \forall t \in \mathcal{T},
	 \label{eqn_reformulation}
\end{align}
where $N_l$ denotes the neuron number in the $l$-th hidden layer. Obviously, the binary variable number introduced by the (\ref{eqn_reformulation}) is the same as the neuron number of the Vio-MLP.

\subsubsection{Loss-MLP}
According to (\ref{eqn_balance}), calculating the net power at the root node, $G_t^\text{root}$, requires the total power loss $p_t^\text{loss}$. However, the value of $p_t^\text{loss}$ is calculated based on the power flow model (\ref{eqn_distflow}), which still requires network parameters. To address this issue, we train another MLP (recorded as ``Loss-MLP") to predict $p_t^\text{loss}$. After training, the Loss-MLP can also be reformulated into a mixed-integer linear form. Since the forward propagation of the Loss-MLP is almost the same as that of the Vio-MLP, we record it as follows:
\begin{align}
\begin{cases}
&\{\text{Eqs. (\ref{eqn_input}) and (\ref{eqn_reformulation})}\}^\text{loss}, \\
&p_t^\text{loss} = (\bm w^{|\mathcal{L}^\text{loss}|+1})^\intercal \bm s_t^{|\mathcal{L}^\text{loss}|} + b^{|\mathcal{L}^\text{loss}|+1},\ \forall t \in \mathcal{T}.
\end{cases} \label{eqn_output_loss}
\end{align}
Here the superscript ``$\text{loss}$" is employed to mark the parameters and variables used in the Loss-MLP. 


Based on the above two MLPs and their mixed integer reformulations, \textbf{P1} can be replicated by the following \textbf{P2}:
\begin{align} 
	&\min_{(\bm p_t^\text{HV}, \bm \lambda_t)_{\forall t \in \mathcal{T}}} \quad \sum_{t \in \mathcal{T}} EC_t \tag{$\textbf{P2}$},\\
	&\begin{array}{r@{\quad}r@{}l@{\quad}l}
		\text{s.t.} &&\text{Eqs.  (\ref{eqn_thermal})-(\ref{eqn_q}), (\ref{eqn_EC})-(\ref{eqn_balance}), (\ref{eqn_x}), (\ref{eqn_input}), and (\ref{eqn_output})-(\ref{eqn_output_loss}).} 
	\end{array} \notag
\end{align}
In \textbf{P2}, Eqs. (\ref{eqn_thermal})-(\ref{eqn_q}) and  (\ref{eqn_x}) describe the relationship between decision variables $(\bm p_t^\text{HV}, \bm \lambda_t)_{\forall t \in \mathcal{T}}$ and input $\bm x_t$; Eqs. (\ref{eqn_EC})-(\ref{eqn_balance}) represent the calculation method of the energy cost $EC_t$; Eqs. (\ref{eqn_input}), (\ref{eqn_output}), and (\ref{eqn_reformulation}) denotes the Vio-MLP; Eq. (\ref{eqn_safe}) defines the safety condition; Eq. (\ref{eqn_output_loss}) represents the Loss-MLP.

Although many published papers have confirmed that constraint learning can accurately replicate power flow constraints without network parameters \cite{9302963,venzke2020neural,9502573,chen2022deep}, this method still faces some challenges. Since neural networks are usually regarded as ``black boxes," it is hard to explain why constraint learning can replicate complex constraints. Meanwhile, it may also be computationally expensive due to the binary variables introduced in (\ref{eqn_reformulation}) and (\ref{eqn_output_loss}). Moreover, its solution may lie in an area where the MLPs have not learned based on historical data, so this solution's reliability may be hard to guarantee. 

\subsection{Piecewise-linearization-based interpretation} \label{sec_PWL}
In this section, we show that constraint learning can be interpreted as a PWL technique for multivariate functions, which explains why it can learn complex constraints. Specifically, we first introduce the concept ``activation pattern." It is a set that contains the activation states of all neurons in an MLP:
\begin{align}
\{\bm o^l, \forall l \in \mathcal{L} \}, \label{eqn_partten}
\end{align}
where vector $\bm o^l$ represents the activation states of all neurons in the $l$-th hidden layer. Its element is defined as:
\begin{enumerate}
\item If $o^l_n=1$, then the $n$-th neuron in hidden layer $l$ is active, and its ReLU's input is non-negative.
\item If $o^l_n=-1$, then the $n$-th neuron in hidden layer $l$ is inactive, and its ReLU's input is negative.
\end{enumerate}
Once the activation pattern is given,  the positive/negative conditions of all the ReLU's inputs are also determined. Since these positive/negative conditions can be regarded as constraints, an activation pattern also bounds a region $\mathcal{R}$:
\begin{align}
\mathcal{R} = \left\{\overbrace{\text{diag}(\bm o^l)}^\text{Activation states} \underbrace{\left(\widehat{\bm W}^l \bm x_t + \widehat{\bm b}^l \right)}_\text{Input of the ReLU} \geq \bm 0, \forall l \in \mathcal{L} \right\}, \label{eqn_activationRegion}
\end{align}
where $\text{diag}(\bm o^l)$ denotes the diagonal matrix formed by the vector $\bm o^l$; $\widehat{\bm W}^l$ and $\widehat{\bm b}^l$ are calculated based on the learning parameters ${\bm W}^l$ and ${\bm b}^l$ in (\ref{eqn_input})-(\ref{eqn_output}):
\begin{align}
\begin{cases}
&\widehat{\bm W}^{1} = \bm W^1, \quad \widehat{\bm b}^{1} = \bm b^1, \\
&\widehat{\bm W}^{l} = \bm W^l \left(\text{diag}(\bm o^{l-1}) \widehat{\bm W}^{l-1}\right), \quad \forall l \in \mathcal{L}/\{1\},\\
&\widehat{\bm b}^{l} = \bm W^l \left(\text{diag}(\bm o^{l-1}) \widehat{\bm b}^{l-1}\right) + {\bm b}^{l}, \quad \forall l \in \mathcal{L}/\{1\}.
\end{cases}
\end{align}
We call this $\mathcal{R}$ as the``activation region" of the given activation pattern. Note that these activation regions are disjoint \cite{sudjianto2020unwrapping}. Then, in every activation region, the output $h_t$ can be calculated based on a linear mapping of input $\bm x_t$:
\begin{align}
h_t = (\widehat{\bm w}^{|\mathcal{L}|+1})^\intercal \bm x_t + \widehat{ b}^{|\mathcal{L}|+1}, \label{eqn_linearMap}
\end{align}
where parameters $\widehat{\bm W}^{|\mathcal{L}|+1}$ and $\widehat{\bm b}^{|\mathcal{L}|+1}$ are defined as:
\begin{align}
\begin{cases}
\widehat{\bm w}^{|\mathcal{L}|+1} = \left( (\bm w^{|\mathcal{L}|+1})^\intercal \left(\text{diag}(\bm o^{|\mathcal{L}|}) \widehat{\bm W}^{|\mathcal{L}|}\right) \right)^{\intercal}, \\
\widehat{b}^{|\mathcal{L}|+1} = \left( (\bm w^{|\mathcal{L}|+1})^\intercal \left(\text{diag}(\bm o^{|\mathcal{L}|}) \widehat{\bm b}^{{|\mathcal{L}|}}\right)\right)^{\intercal} + {b}^{|\mathcal{L}|+1}.
\end{cases} \label{eqn_parameter_output}
\end{align}
Note that each neuron has two activation states. Thus, one MLP can have $|\mathcal{K}|=2^{N^\text{Total}}$ activation regions, where $N^\text{Total}$ is the total neuron number of this MLP.

Eq. (\ref{eqn_safe}) requires that the measurement of constraint violations $h_t$ should be non-positive. If we use $k \in \mathcal{K}$ to index activation regions, then the feasible region in each activation region, i.e., $\mathcal{F}_k$, can be expressed as:
\begin{align}
\mathcal{F}_k = \left\{ \bm x_t \left| \begin{matrix}  &\bm x_t \in \mathcal{R}_{k}, \\ 
&h_{k,t} = (\widehat{\bm w}_k^{|\mathcal{L}|+1})^\intercal \bm x_t + \widehat{ b}_k^{|\mathcal{L}|+1} \leq 0, \end{matrix} \right. \right\}, \forall k \in \mathcal{K}. \label{eqn_f}
\end{align}
Since all constraints in (\ref{eqn_f}) is linear, set $\mathcal{F}_k$ can be further expressed as a polytope:
\begin{align}
\mathcal{F}_k = \left\{ \bm x_t \left| \bm A_k \bm x_t \leq \bm \beta_k \right.  \right\}, \forall k \in \mathcal{K}, \label{eqn_f2}
\end{align}
where $\bm A_k$ and $\bm \beta_k$ are the polytope's parameters defined based on (\ref{eqn_f}). Then, the feasible region formed by the conventional constraint learning, i.e., (\ref{eqn_input})-(\ref{eqn_safe}), can be rewritten as a union of $\mathcal{F}_k$ in all activation regions:
\begin{align}
\bm x_t \in \cup_{k \in \mathcal{K}} \mathcal{F}_k. \label{eqn_union}
\end{align}

\begin{figure}
		\vspace{-4mm}
	\centering
	{\includegraphics[width=0.98\columnwidth]{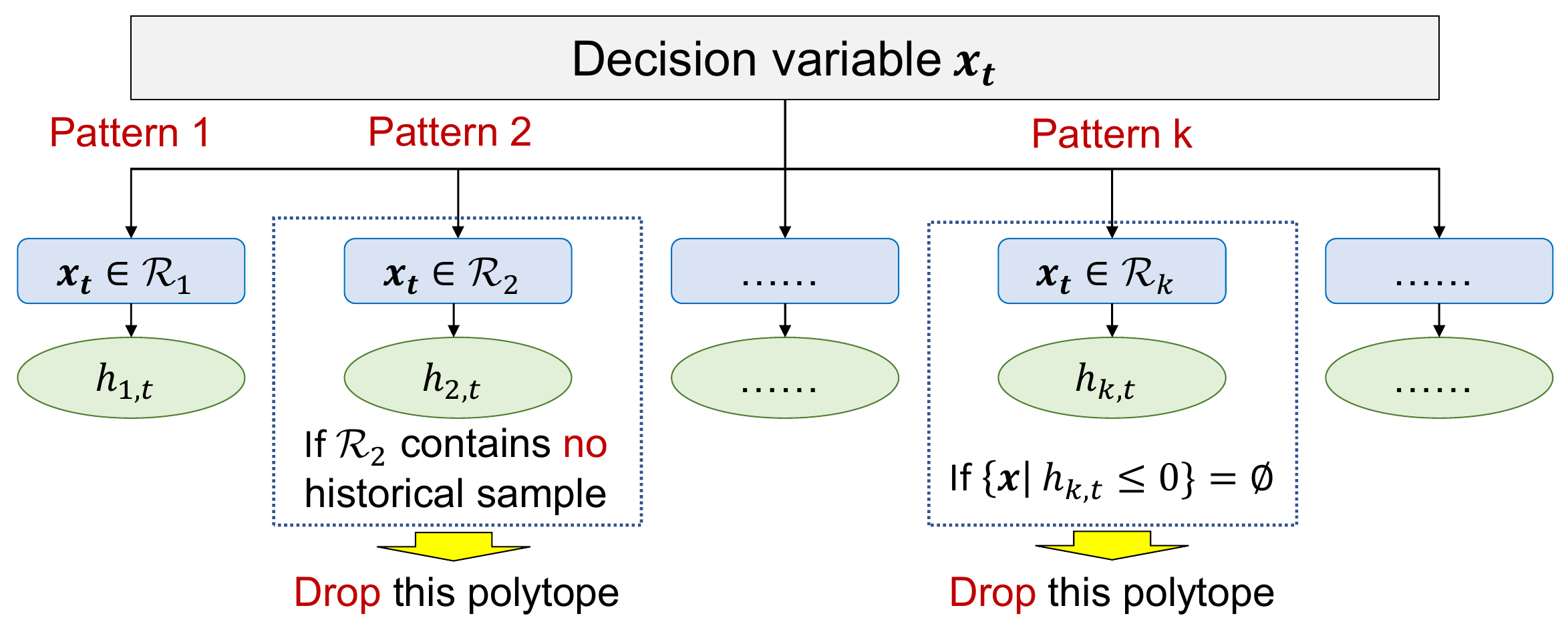}} 
	\vspace{-2mm}
 	\caption{PWL-based interpretation of constraint learning. For a given input $\bm x$, we first check that it belongs to which activation region $\mathcal{R}_k$ according to its activation pattern. Then, the output $h_t$ can be calculated based on the corresponding linear mapping.}
	\label{fig_tree}
	\vspace{-4mm}
\end{figure}

\begin{remark} \label{remark_2}
According to (\ref{eqn_partten})-(\ref{eqn_union}), the process of constraint learning can be summarized by Fig. \ref{fig_tree}: Based on different activation patterns, the trained MLP first partitions its input domain into disjoint linear pieces (activation regions). In each piece, the output of the MLP is linearly dependent on the input. Then, the feasible region of the original problem can be approximated by a union of polytopes. This manner is very similar to the PWL for multivariate functions.
\end{remark}

\begin{remark}
Increasing the number of linear pieces (i.e., activation regions) can reduce the approximation errors of the PWL (i.e., constraint learning). Meanwhile, reference \cite{hanin2019deep} pointed out that the activation region number of an MLP is positively correlated with its neuron number. Thus, constraint learning can accurately replicate constraints once the MLP has a proper number of neurons.
\end{remark}

\begin{figure}
		\vspace{-4mm}
	\subfigbottomskip=-4pt
	\subfigcapskip=-4pt
	\centering
	\subfigure[]{\includegraphics[width=0.49\columnwidth]{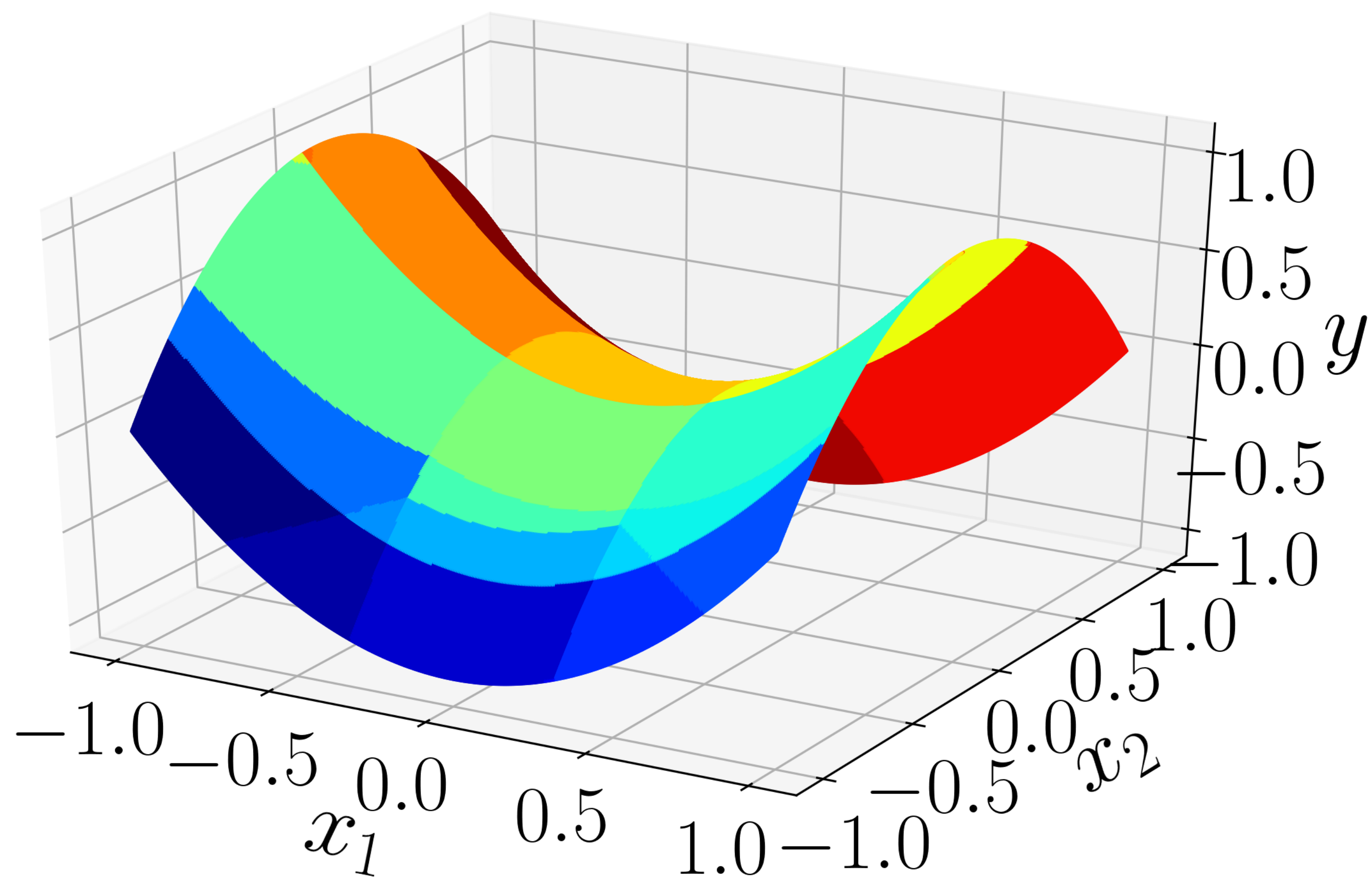}} 
	\subfigure[]{\includegraphics[width=0.49\columnwidth]{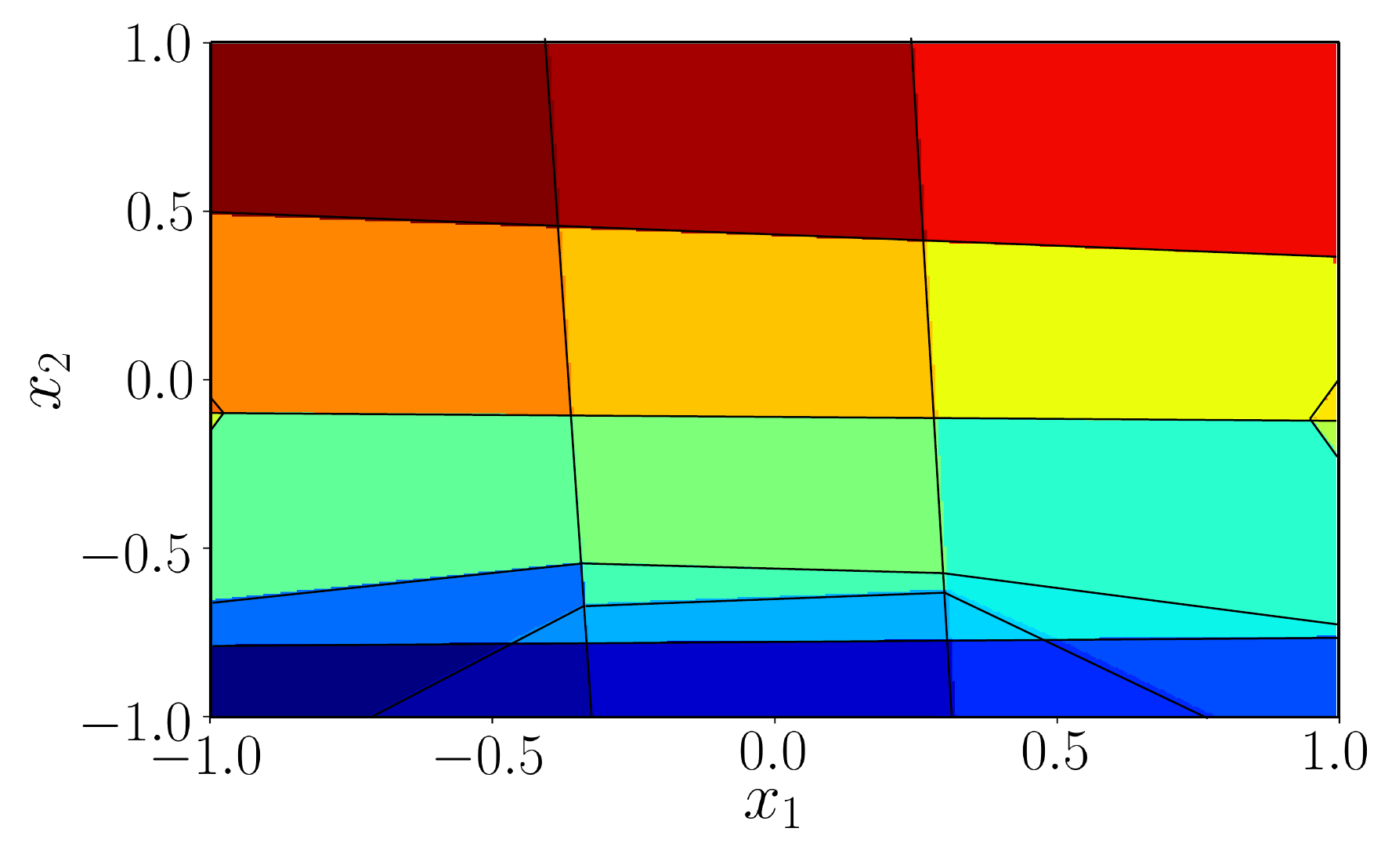}}
	\subfigure[]{\includegraphics[width=0.49\columnwidth]{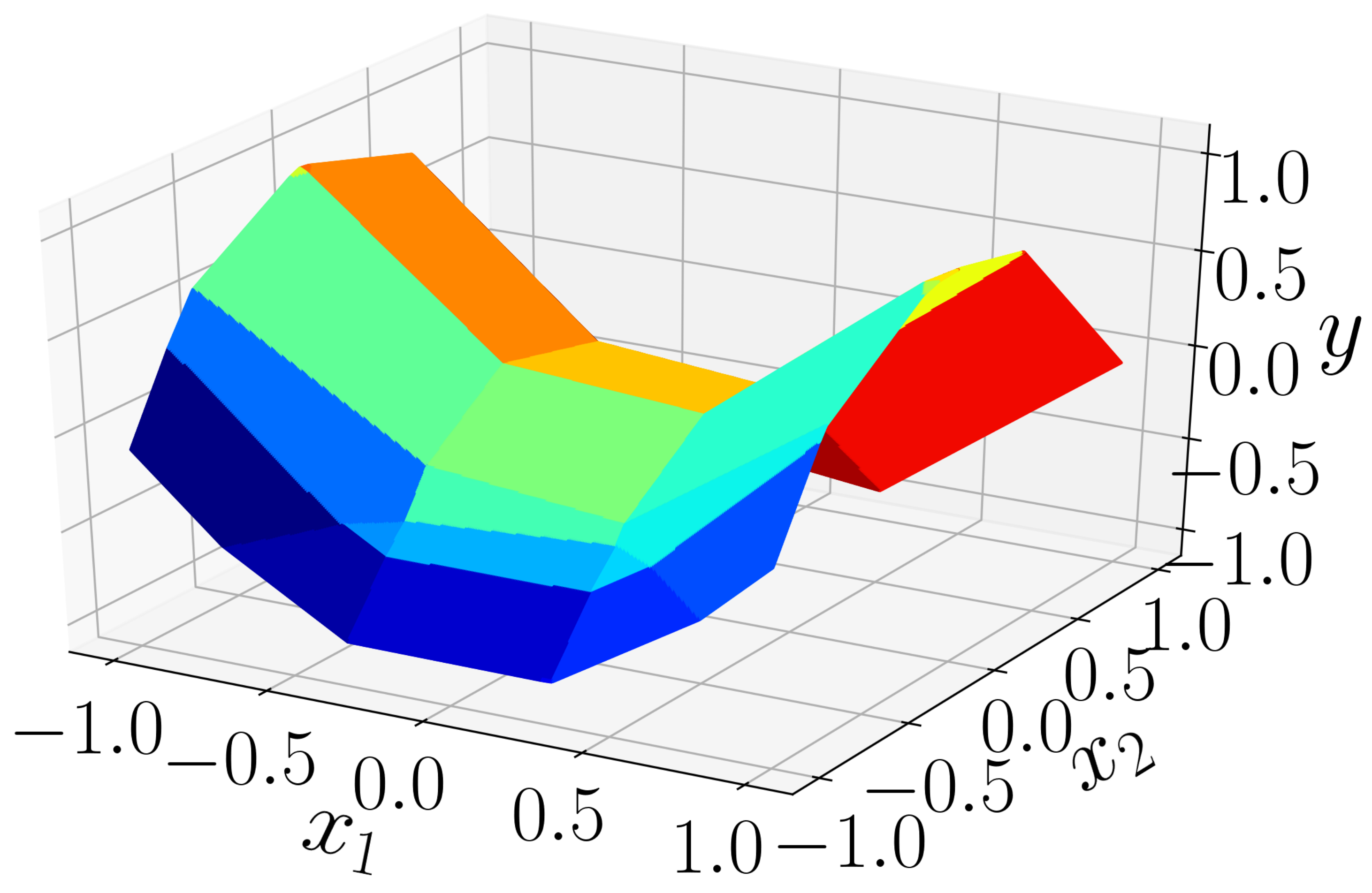}}
	\subfigure[]{\includegraphics[width=0.49\columnwidth]{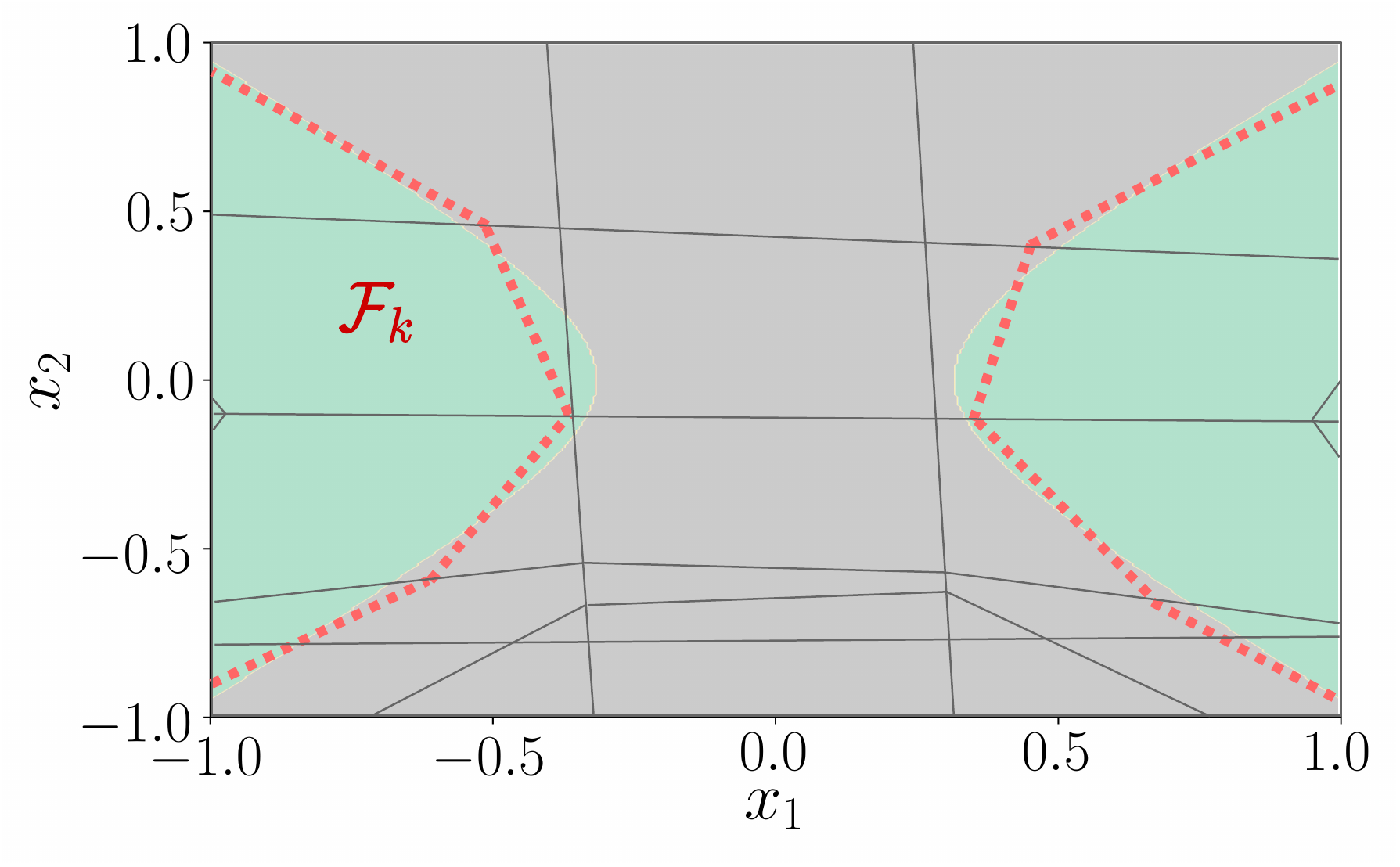}}
 	\caption{An example to explain the PWL-based interpretation, where a MLP with 15 neurons are trained to replicate the constraint $y=x_1^2-x_2^2 \geq 0.1$. 
 	(a) Function graph of $y=x_1^2-x_2^2$, (b) different activation regions (distinguished by colors), (c) output of the trained MLP,  and (d) the original feasible region bounded by $y \geq 0.1$ (green area) and its approximation provided by constraint learning (areas bounded by the red dash lines and black lines). The red dash lines in (d) represent the positions where the predicted $y$ equals to 0.1.}
	\label{fig_activationRegion}
	\vspace{-4mm}
\end{figure}

We provide an example in Fig. \ref{fig_activationRegion} to further explain this PWL-based interpretation. In this example, constraint learning trains an MLP with three hidden layers to replicate a non-convex constraint $y=x_1^2-x_2^2 \geq 0.1$ with $x_1, x_2 \in [-1, 1]$. 
Each hidden layer contains five neurons. Fig. \ref{fig_activationRegion}(a) shows the graph of the function $y=x_1^2-x_2^2$. Fig. \ref{fig_activationRegion}(b) illustrates different activation regions of the trained MLP, which are distinguished by colors. Fig. \ref{fig_activationRegion}(c) demonstrates the output of the trained MLP. The output surface is composed of multiple planes. These planes demonstrate the linear mappings from $x_1$ and $x_2$ to $y$. Fig. \ref{fig_activationRegion}(d) compares the actual feasible region $y \geq 0.1$ (green areas) and the approximation provided by constraint learning (areas bounded by the red dash and black solid lines). Obviously, this feasible region can be accurately represented even if the MLP only has 15 neurons.

\subsection{Problem simplification} \label{sec_prune}
\textbf{P2} is a mixed-integer linear program, so the Branch-and-Bound algorithm can solve it with guaranteed optimality. In this algorithm, the linear problems relaxed from \textbf{P2} are recursively solved to update the lower and upper bounds of the solution. If we use $N^\text{Total,Vio}$ and $N^\text{Total,Loss}$ to represent the total neuron numbers of the Vio-MLP and Loss-MLP, then \textbf{P2} totally introduces $N^\text{Total,Vio}+N^\text{Total,Loss}$ binary variables. Hence, we need to solve $2^{N^\text{Total,Vio}+N^\text{Total,Loss}}$ linear problems in the worst case, leading to a huge computational burden.


We propose a two-step simplification method to address this problem. The first step leverages existing constraints and historical samples to prune away unnecessary polytopes in the union (\ref{eqn_union}). The second step removes redundant linear constraints for each retained polytope. Then, the computational efficiency can be enhanced. Moreover, since this simplification method can avoid extrapolating significantly from the historical samples, the reliability of solutions can also be enhanced. 

\subsubsection{Pruning away unnecessary polytopes}

In fact, many activation regions are unnecessary. On the one hand, most activation regions do not cover any historical sample. In these regions, the MLP may not guarantee prediction accuracy. Thus, if the solution of constraint learning lies in these regions, it may be unreliable. On the other hand, some activation regions may not contain any feasible solution, e.g.,  $\mathcal{F}_k = \varnothing$. As a result, it is preferable to remove these unnecessary activation regions in advance, as shown in Fig. \ref{fig_tree}. Then, the computational burden can be reduced. For example, the MLP in Fig. \ref{fig_activationRegion} has 15 neurons, so it totally has $2^{15}=32718$ activation regions. Thus, the branch-and-bound may solve $32718$ relaxed linear problems in the worst case. Nevertheless, only $14$ activation regions cover historical samples and feasible solutions simultaneously. Hence, we only need to solve at most $14$ relaxed linear problems after pruning.

To identify all the unnecessary activation regions in our problem, we first check that each historical sample belongs to which activation region. Then, the activation regions that cover historical samples can be recognized. Next, the following optimization problem \textbf{P-S1} is solved to judge whether an activation region contains feasible solutions or not:
\begin{align} 
	&\min_{\bm x_t \in \mathcal{X}} \  0, \quad \text{s.t. }  \bm x_t \in \mathcal{F}_k \tag{$\textbf{P-S1}$},
\end{align}
where $\mathcal{X}$ is the domain of nodal power injections $\bm x_t$, which can be approximated based on the historical data. If \textbf{P-S1} has solutions, then the $k$-th activation region contains feasible solutions of \textbf{P2}. Otherwise, it does not contain any feasible solution and can be removed. Once all legal activation regions are found, constraint (\ref{eqn_union}) can be replaced by:
\begin{align}
\bm x_t \in \cup_{k \in \mathcal{K}^\text{Fe}} \mathcal{F}_k, \label{eqn_union_prune}
\end{align}
where $\mathcal{K}^\text{Fe}$ is the index set of the legal activation regions. 
Based on the Big-M method, constraint (\ref{eqn_union_prune}) can be further reformulated into a mixed-integer linear form:
\begin{align}
&\bm A_k \bm x_{k,t} \leq z_{k,t} \beta_k, \ \forall k \in \mathcal{K}^\text{Fe}, \label{eqn_MIP_prune} \\
&-M z_{k,t} \leq \bm x_{k,t} \leq M z_{k,t}, \ \forall k \in \mathcal{K}^\text{Fe}, \label{eqn_MIP_prune2} \\
& \sum_{k \in \mathcal{K}^\text{Fe}} x_{k,t} = \bm x_t,\  \sum_{k \in \mathcal{K}^\text{Fe}} z_{k,t} = 1,\  \bm z_t \in \{0,1\}^{|\mathcal{K}^\text{Fe}|}. \label{eqn_binary_prune}
\end{align}
Obviously, due to the constraint $\sum_{k \in \mathcal{K}^\text{Fe}} z_{k,t} = 1$, at most $|\mathcal{K}^\text{Fe}|$ linear problems are introduced in the worst case, which is much smaller compared to the conventional constraint learning \textbf{P2}.\footnote{
State-of-art solvers (e.g., GUROBI) usually have a pre-solving function. This function can leverage existing constraints to remove the searching space containing no feasible solution in advance. The proposed simplification method not only removes the space without any feasible solution but also drops the place that covers no historical sample. Therefore, it can remove more space compared to the pre-solving function. As a result, it has the potential to achieve better computational efficiency. 
}
Thus, 

\subsubsection{Removing redundant constraints}
According to (\ref{eqn_f}), each polytope $\mathcal{F}_k$ is bounded by multiple linear constraints. Since (\ref{eqn_MIP_prune}) involve multiple polytopes, it may introduce large-scale constraints and deteriorate the computational performance. To address this issue, we further design a step to remove redundant constraints for each polytope $\mathcal{F}_k = \left\{ \bm x_t \left| \bm A_k \bm x_t \leq \bm \beta_k \right.  \right\}$. Specifically, for a given polytope $\mathcal{F}_k$, we can solve the following problem to identify whether its $i$-th constraint is redundant or not:
\begin{align} 
	\rho_k^{(i)}=&\max_{\bm x_t \in \mathcal{X}} \  \bm A_k^{(i)} \bm x_t, \tag{$\textbf{P-S2}$}\\
	&\begin{array}{r@{\quad}r@{}l@{\quad}l}
		\text{s.t.} &&\bm A_k^{(m)} \bm x_t \leq  \beta_k^{(m)}, \forall m \in \mathcal{M}/\{i\}, 
	\end{array} \notag
\end{align}
where $\bm A_k^{(m)}$ and $\beta_k^{(m)}$ represent the $m$-th rows of $\bm A_k$ and $\bm \beta_k$; $m \in \mathcal{M}$ is the index of constraints in $\mathcal{F}_k$. If the solution $\rho_k^{(i)} \leq {\beta}_k^{(i)}$, then the $i$-th linear constraint, i.e., $\bm A_k^{(i)} \bm x_t \leq \beta_k^{(i)}$, is redundant. Otherwise, it may be active.

\begin{figure}
		\vspace{-4mm}
	\centering
	{\includegraphics[width=0.8\columnwidth]{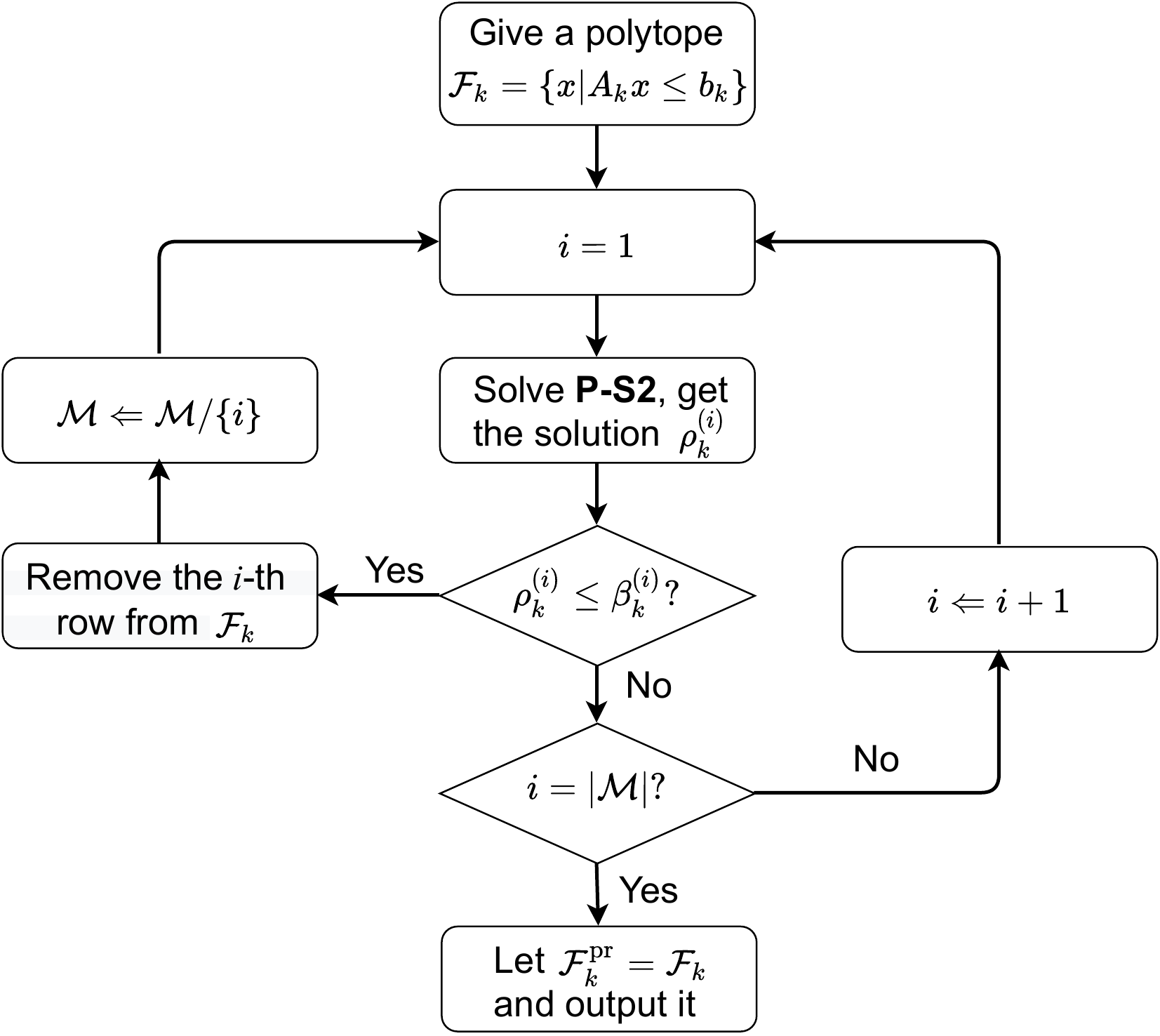}} 
	\vspace{-2mm}
 	\caption{Algorithm to remove redundant constraints in $\mathcal{F}_{k}$. The output is recorded as $\mathcal{F}^{\text{pr}}_k=\{x|\bm A^{\text{pr}}_k \bm x_t \leq \bm b^{\text{pr}}_k\}$.}
	\label{fig_redundant}
	\vspace{-4mm}
\end{figure}

It should be noted that removing one constraint may affect the redundancy of another constraint. For example, we may have $\rho_k^{(i)} \leq {\beta}_k^{(i)}$ and $\rho_k^{(j)} \leq {\beta}_k^{(j)}$ simultaneously for a given $\mathcal{F}_k$. It seems that both the $i$-th and $j$-th constraints are redundant. However, once the $i$-th constraint is removed, the $j$-th constraint may not be redundant anymore in the new polytope $\mathcal{F}_k' = \{\bm A_k^{(m)} \bm x_t \leq \beta_k^{(m)}, \forall m \in \mathcal{M}/\{i\}\}$. Based on this observation, we design a heuristic algorithm to identify all redundant constraints, as shown in Fig. \ref{fig_redundant}.
Its detailed steps can be summarized as: 
\begin{enumerate}
\item For a given polytope, solve \textbf{P-S2} with $i \in \mathcal{M}$ in turn.
\item Once the first redundant constraint is found (indexed by $j$), remove it from $\mathcal{F}^k$, resulting in a new polytope. 
\item Repeat Steps 1 and 2 until no redundant constraint can be found. Output the newest polytope (recorded as $\mathcal{F}^{\text{pr}}_k=\{x|\bm A^{\text{pr}}_k\bm x\leq \bm b^{\text{pr}}_k\}$).
\end{enumerate}
Then, Eq. (\ref{eqn_MIP_prune}) can be equivalently simplified as:
\begin{align}
\bm A^{\text{pr}}_k\bm x\leq z_{k,t} \bm b^{\text{pr}}_k, \  \forall k \in \mathcal{K}^\text{Fe}. \label{eqn_redundant_prune}
\end{align}

\subsection{Summary of the efficient constraint learning}
\subsubsection{Final formulation}
After applying the proposed two-step pruning, \textbf{P2} can be replaced by:
\begin{align} 
	&\min_{(\bm p_t^\text{HV}, \bm \lambda_t)_{\forall t \in \mathcal{T}}} \quad \sum_{t \in \mathcal{T}} EC_t \tag{$\textbf{P3}$},\\
	&\begin{array}{r@{\quad}r@{}l@{\quad}l}
		\text{s.t.} &&\text{Eqs.  (\ref{eqn_thermal})-(\ref{eqn_q}), (\ref{eqn_EC})-(\ref{eqn_balance}), (\ref{eqn_output_loss}), \{(\ref{eqn_MIP_prune2})-(\ref{eqn_redundant_prune})\}$^\text{Vio}$.} 
	\end{array} \notag
\end{align}
where symbol \{(\ref{eqn_MIP_prune2})-(\ref{eqn_redundant_prune})\}$^\text{Vio}$ represents that the Vio-MLP is replicated by the efficient reformulation (\ref{eqn_MIP_prune2})-(\ref{eqn_redundant_prune}).\footnote{Note we do not apply the simplification method to the Loss-MLP. The power loss $p_t^\text{loss}$ is only used to estimate the net power $G_t^\text{root}$ based on (\ref{eqn_balance}). Since $p_t^\text{loss}$ is usually much smaller than the summation of nodal active power injections, the estimated $G_t^\text{root}$ can still be accurate enough even if the accuracy of predicting $p_t^\text{loss}$ is not very high. In other words, a small number of neurons is already enough for the Loss-MLP. Thus, it is unnecessary to further simplify the Loss-MLP because its structure already is simple enough.} Compared to \textbf{P2}, \textbf{P3} introduces much fewer activation regions. Moreover, the redundant linear constraints are also dropped. Hence, its computational efficiency is significantly enhanced. Meanwhile, since the simplification method prevents excessive exploration, \textbf{P3}'s feasible region does not include those activation regions that the MLP has not learned. Hence, the reliability of the solutions is also improved.

\subsubsection{Measurement requirement}
The proposed method relies on data measurements from distribution networks. According to (\ref{eqn_x})-(\ref{eqn_h}), the required historical samples include nodal power injections, bus voltages, and branch power flows. In practice, power injections are usually monitored because they determine the electricity bill of users. However, since the measurement redundancy on a distribution system is usually low, not only some network parameters (i.e., topology and impedance) but also parts of bus voltage and branch power flows are not measured. Nevertheless, we can still calculate the measurement of constraint violations $h_t$ with those measured data based on (\ref{eqn_h}). Then, the proposed method can derive a strategy that at least guarantees the security of the measured buses and branches. On the contrary, other model-based methods or data-driven models (e.g., optimize-then-learning methods) can hardly derive a solution without sufficient network parameters. Hence, the proposed model is still meaningful for distribution systems with low measurement redundancy.

\section{Case study} \label{sec_case}
\subsection{Simulation set up}
Two different case studies are implemented based on the IEEE 33- and 123-bus systems to verify the effectiveness of the proposed efficient constraint learning. In both cases, the time interval and optimization horizon are set as one hour and 24 hours, respectively. The parameters of buildings and HVAC systems are summarized in Table \ref{tab_parameter}.  
Other parameters will be introduced in Sections \ref{sec_33Bus} and \ref{sec_123Bus}.

\begin{table}
	\small
	\centering
	\caption{Parameters of buildings and HVAC systems}
	\vspace{-2mm}
	\begin{tabular}{cccc}
		\hline
		\rule{0pt}{11pt}		
		Parameters & Value &Parameters & Value\\
		\hline
		\rule{0pt}{10pt}
		$C_i$ &  1MWh/{\textcelsius} & $\bm \theta^\text{min}$& 24{\textcelsius}\\
		$\phi_i$ & 0.98 & $\bm \theta^\text{max}$& 28{\textcelsius} \\
		$\text{COP}_i$ &  6 &  $\bm p_{t}^\text{HV,max}$ & 0.1MW \\
		\hline
	\end{tabular}\label{tab_parameter}
	\vspace{-4mm}
\end{table}

The historical operational data of ADNs are simulated based on Pandapower, a power system simulation toolbox in Python \cite{8344496}. The Pandapower is based on the full AC power flow model. Specifically, we first generate 20,000 samples of the input $\bm x_t$. By giving these samples to Pandapower, bus voltages and power flows can be calculated. Based on (\ref{eqn_loss}) and (\ref{eqn_h}), the power loss $p_t^\text{loss}$ and measurement of the constraint violation $h_t$ can be obtained. We gather all generated samples in the form of ($\bm x_t$, $h_t$) and ($\bm x_t$, $p_t^\text{loss}$) to construct training sets for the Vio-MLP and Loss-MLP, respectively. All these samples have been uploaded in \cite{samples2022}.

All numerical experiments are implemented on an Intel(R) 8700 3.20GHz CPU with 16 GB memory. Our MLPs are implemented and trained by Pytorch. 
We employ CVXPY to build optimization problems and GUROBI to solve them.

\subsection{Benchmarks}
To demonstrate the benefits of the proposed method, two benchmarks are introduced:
\begin{enumerate}
	\item \textbf{B1}: Original constraint learning used in \cite{9502573,chen2022deep};
	\item \textbf{B2}: SOCP relaxation of the DistFlow used in \cite{9372881,9112232,9628048}.
\end{enumerate}
Note that the network parameters are unknown in the proposed method and \textbf{B1}, while it is known in \textbf{B2}. Meanwhile, the MLPs used in \textbf{B1} are the same as those in the proposed method. 

\subsection{Case study based on the IEEE 33-bus system} \label{sec_33Bus}
\subsubsection{Parameter setting} 
\begin{figure}
		\vspace{-4mm}
	\centering
	{\includegraphics[width=0.9\columnwidth]{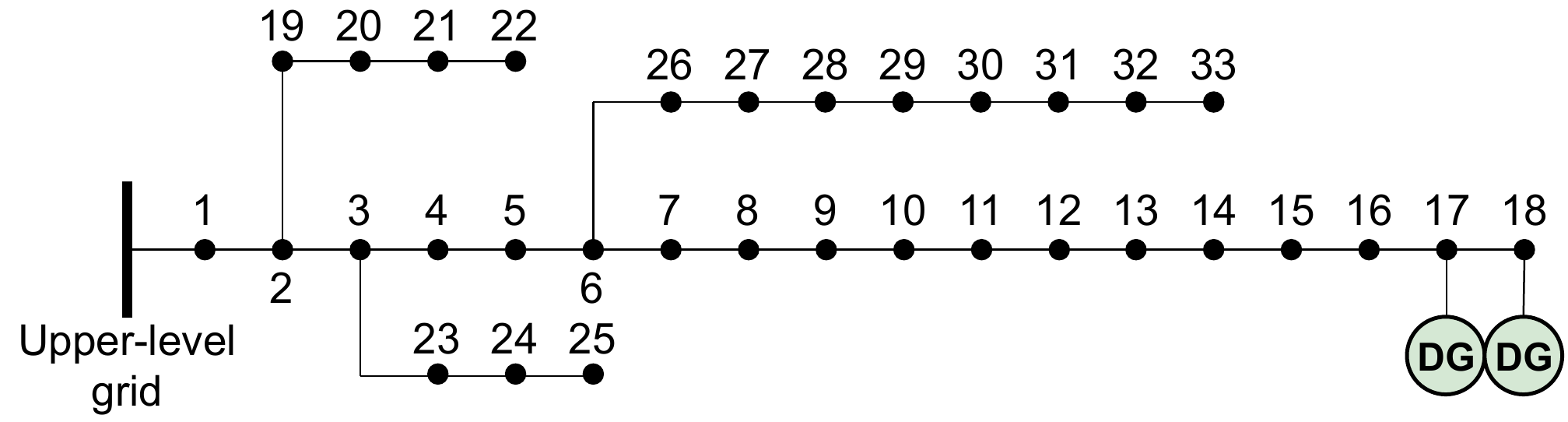}} 
	\vspace{-4mm}
 	\caption{Structure of the IEEE 33-bus system, where ``DG" denotes the distributed renewable generator.}
	\label{fig_33Bus}
\end{figure}

\begin{figure}
		\vspace{-4mm}
	\subfigbottomskip=-4pt
	\subfigcapskip=-4pt
	\centering
	\subfigure[]{\includegraphics[width=0.49\columnwidth]{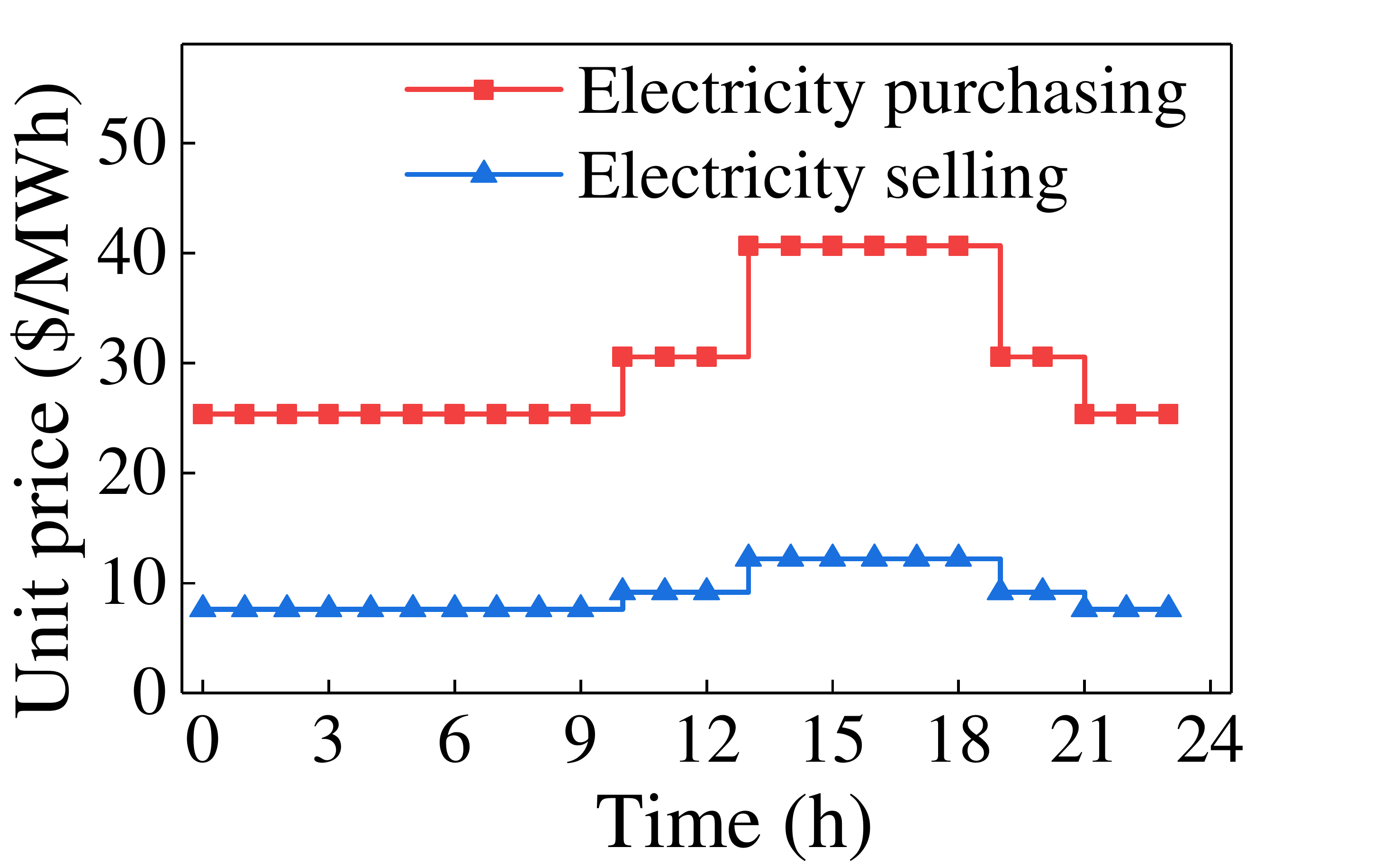}} 
	\subfigure[]{\includegraphics[width=0.49\columnwidth]{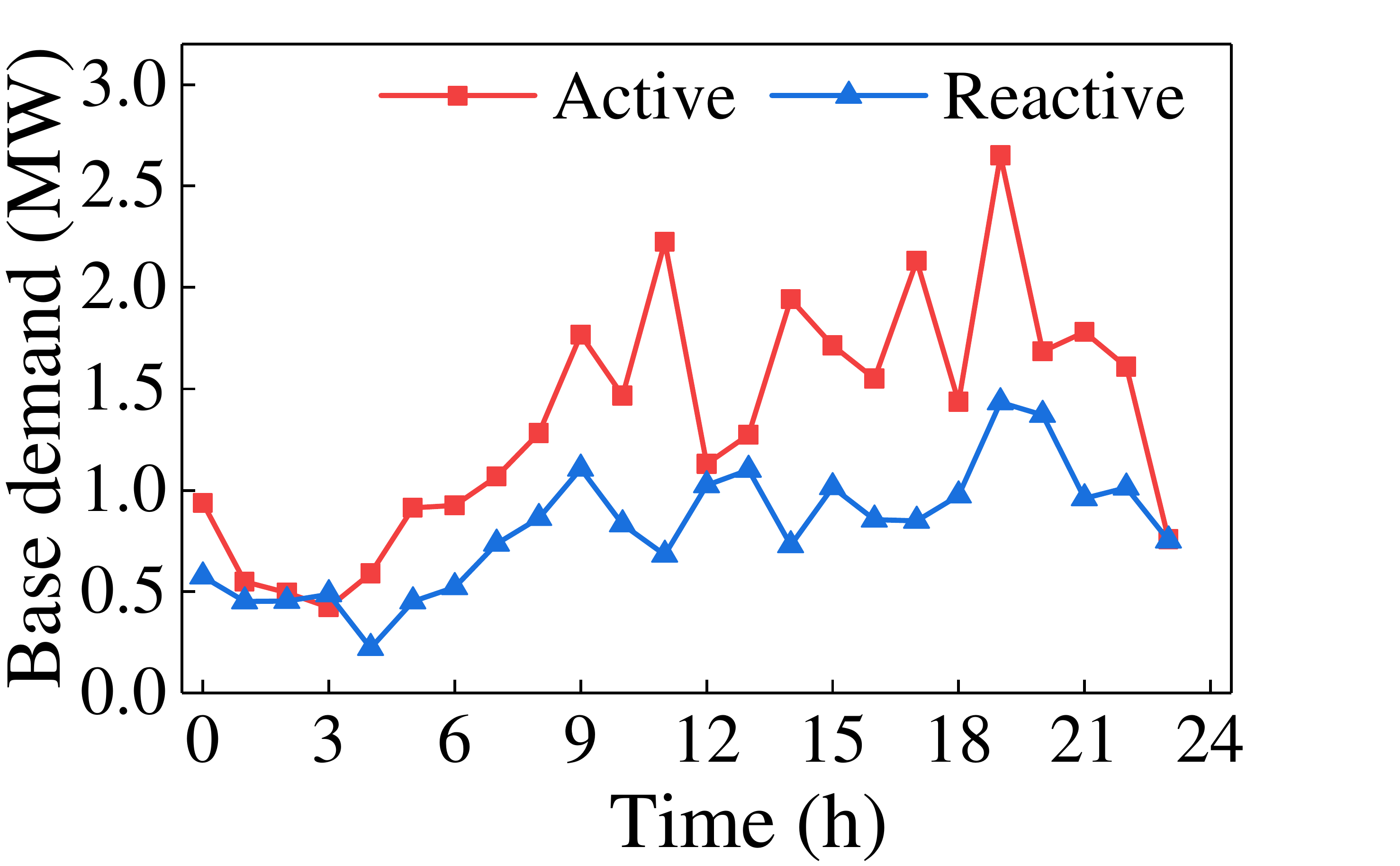}}
	\subfigure[]{\includegraphics[width=0.49\columnwidth]{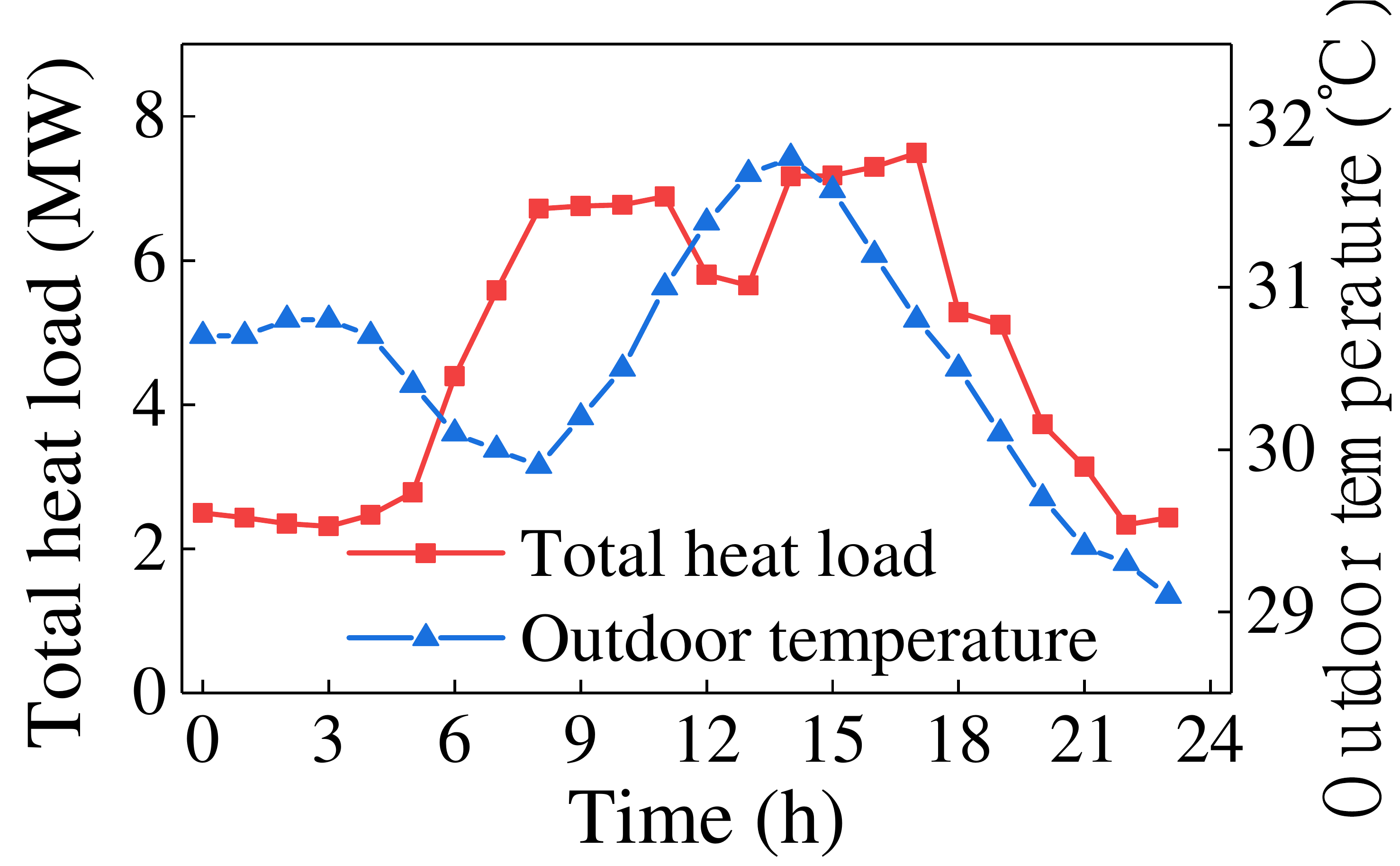}}
	\subfigure[]{\includegraphics[width=0.49\columnwidth]{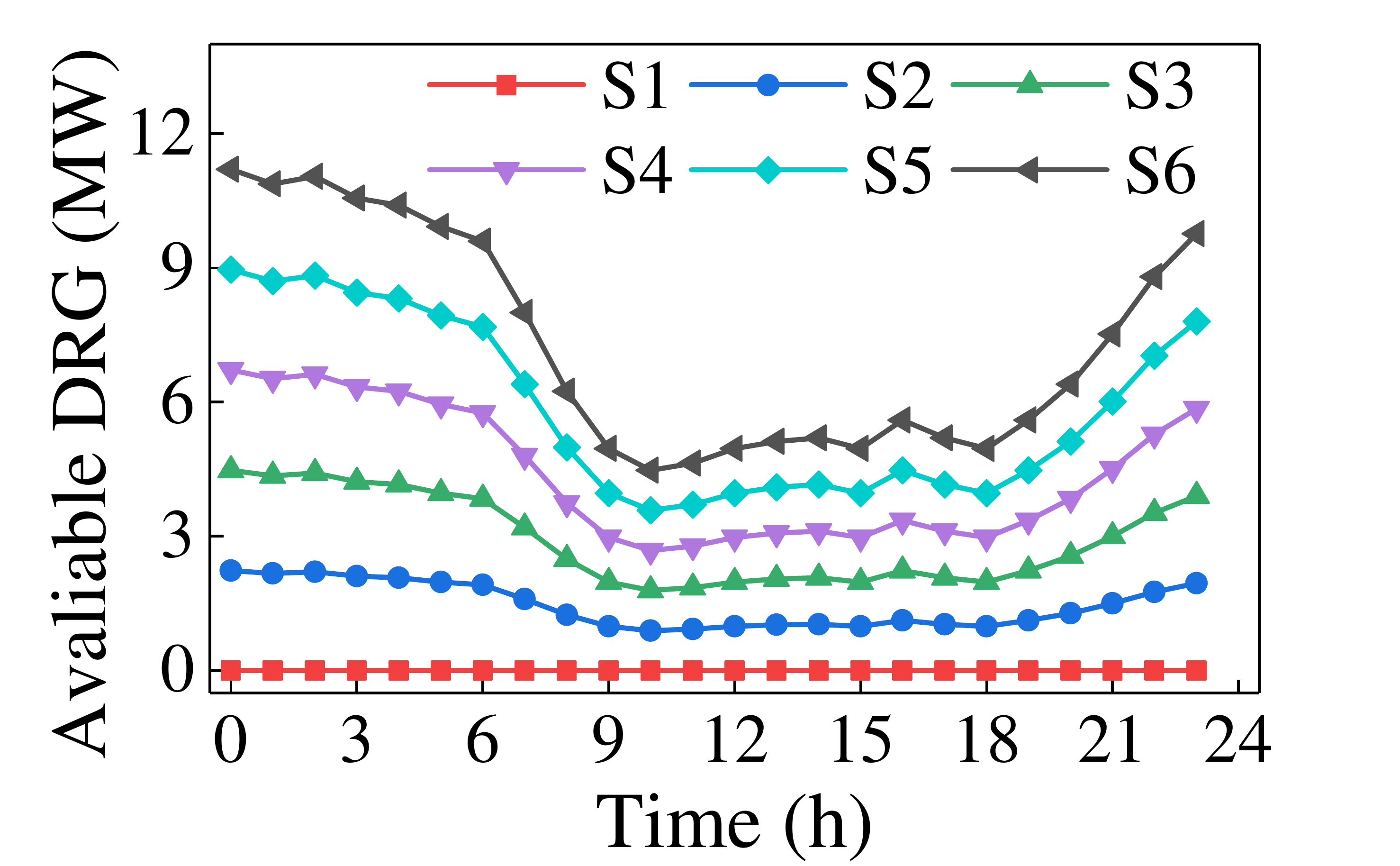}}
 	\caption{Parameters of (a) unit prices of electricity purchasing and selling, (a) total active and reactive base power demands, (c) total heat load contributed by indoor sources, and outdoor temperature, and (d) six scenarios with different available DRG in the case based on the IEEE 33-bus system.}
	\label{fig_case_study_33Bus}
	\vspace{-4mm}
\end{figure}

The structure of the 33-bus system is shown in Fig. \ref{fig_33Bus}. It contains two distributed generators (marked by ``DG"). The voltage of its root node is 12.66kV. The bus voltages are restricted in [0.9p.u., 1.1p.u.], while the maximum allowable value of the line's apparent power flow is 4MW. The heat transfer coefficient of the building, i.e., $g_i$, is 0.03MW/{\textcelsius}. The unit prices for electricity purchasing and selling, base power demands (the total demands except for HVAC loads), total heat loads contributed by indoor sources, and outdoor temperature are illustrated in Figs. \ref{fig_case_study_33Bus}(a)-(c). To better clarify the benefits of the proposed method, we also implement six scenarios with different amounts of available DRG, as shown in Fig. \ref{fig_case_study_33Bus}(d). 

\subsubsection{Performance of MLPs}
\begin{figure}
		\vspace{-4mm}
	\subfigbottomskip=-4pt
	\subfigcapskip=-4pt
	\centering
	\subfigure[]{\includegraphics[width=0.5\columnwidth]{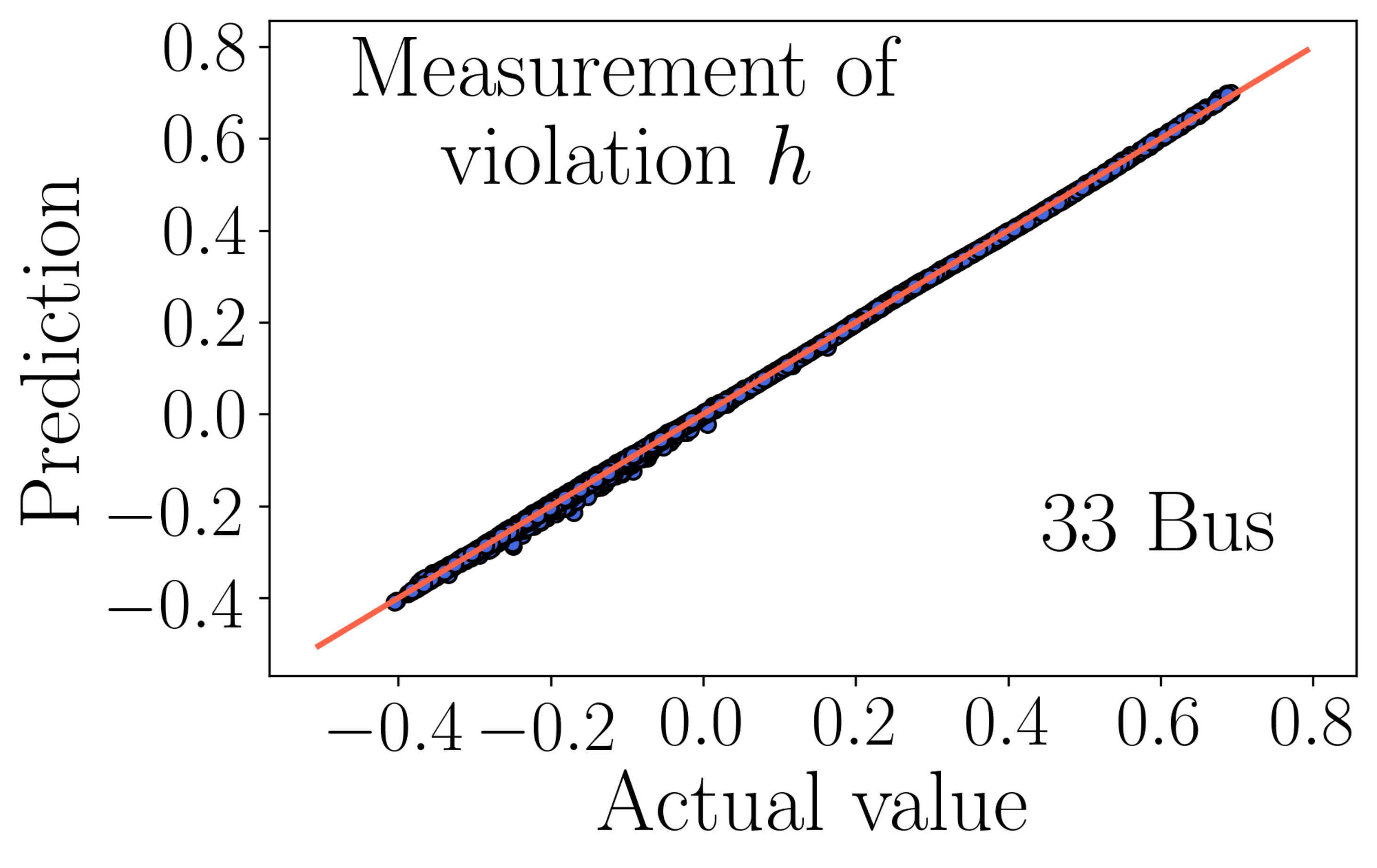}} 
	\subfigure[]{\includegraphics[width=0.48\columnwidth]{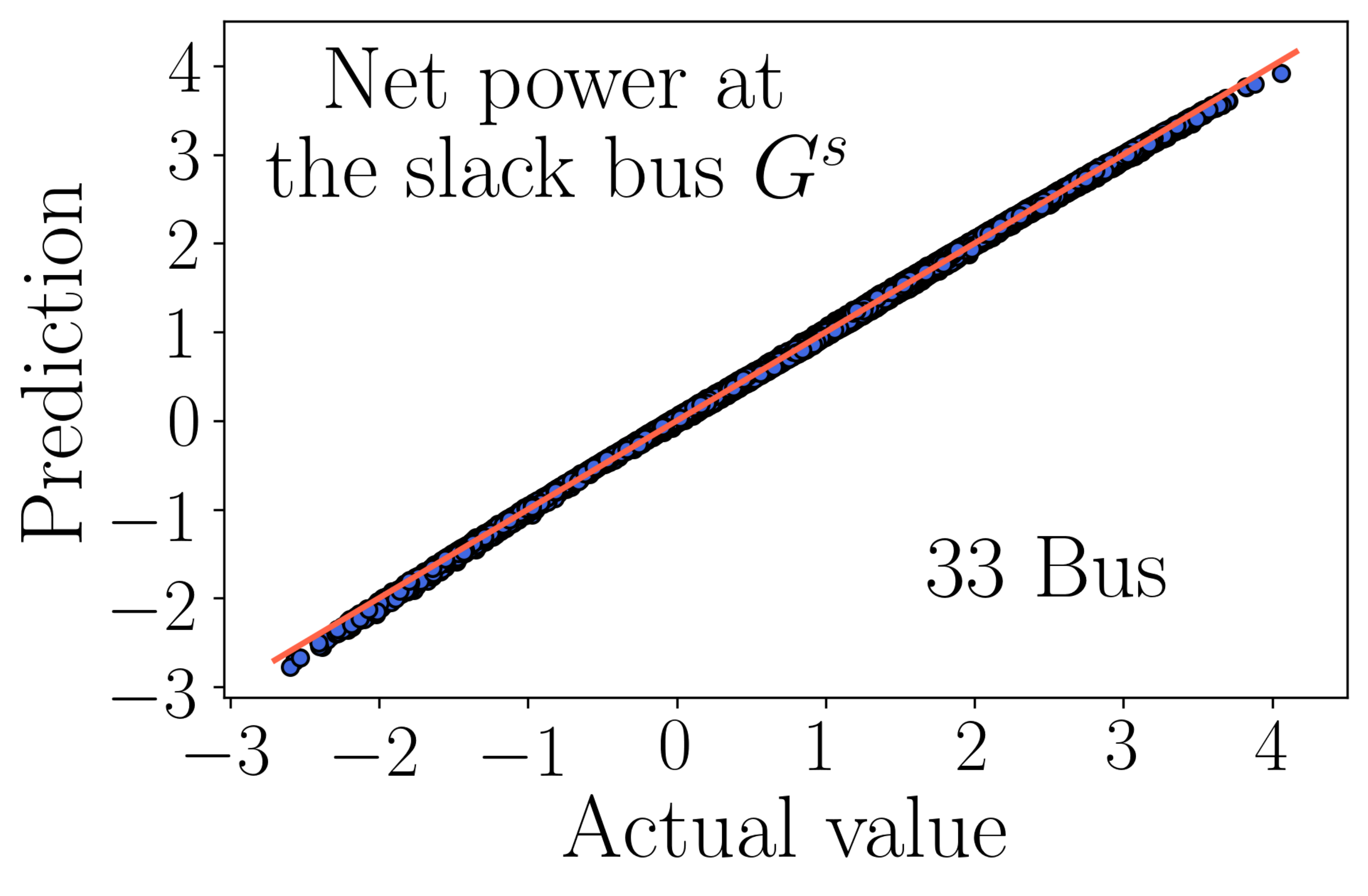}}
	\vspace{-4mm}
 	\caption{The actual and predicted values of (a) the measurement of constraint violation $h_t$ and (b) net power at the root node in the 33-bus case.}
	\label{fig_MLP_performance}
	\vspace{-4mm}
\end{figure}

Fig. \ref{fig_MLP_performance} demonstrates the prediction accuracy of the Vio-MLP and Loss-MLP in the 33-bus case. Note that we directly plot the actual and estimated values of the net power at the root node, i.e., $G_t^\text{root}$, because the power loss is only used to estimate $G_t^\text{root}$. The horizontal and vertical axes represent the actual values and corresponding estimations, while the red line is the position where the estimation equals the actual value. The structure of the Vio-MLP is (6, 6, 6), i.e., three hidden layers with six neurons in each layer, while it is (3, 3, 3) in the Loss-MLP. Obviously, all the samples are pretty close to the corresponding red lines, which indicates the high prediction accuracy of MLPs. 

\subsubsection{Effectiveness of the proposed simplification method}
Table \ref{tab_prune} compares the activation region numbers (column ``Activation regions") and average linear constraint numbers in each activation region (column ``Average constraints") before and after applying the proposed simplification method. Since the Vio-MLP contains 18 neurons, it totally has $2^{18}=262,144$ activation regions. Meanwhile, according to (\ref{eqn_f2}), each activation region is formed by $18+1=19$ linear constraints. The first eighteen constraints form the activation region, i.e., $\bm x_t \in \mathcal{R}_k$, while the last one requires the violation measurement $h_{k,t}$ is non-positive.
After simplification, the above two numbers decrease to 10 and 7.8, respectively. This demonstrates the effectiveness of our simplification method.
\begin{table}
	\small
	\centering
	\caption{Comparison between the cases with and without pruning}
	\vspace{-2mm}
	\begin{threeparttable} 
	\begin{tabular}{ccc}
		\hline
		\rule{0pt}{11pt}		
		 \textbf{Numbers} & \textbf{Activation regions}  & \textbf{Average constraints} \\
		\hline
		\rule{0pt}{10pt}
		Before pruning & 262,144 & 19\\
		After pruning & 10 & 7.8\\
		\hline
	\end{tabular}\label{tab_prune}
	\end{threeparttable}
	\vspace{-4mm}
\end{table}

\subsubsection{Optimality and feasibility}
\begin{figure}
	\subfigbottomskip=-7pt
	\subfigcapskip=-7pt
	\centering
	\subfigure[]{\includegraphics[width=0.9\columnwidth]{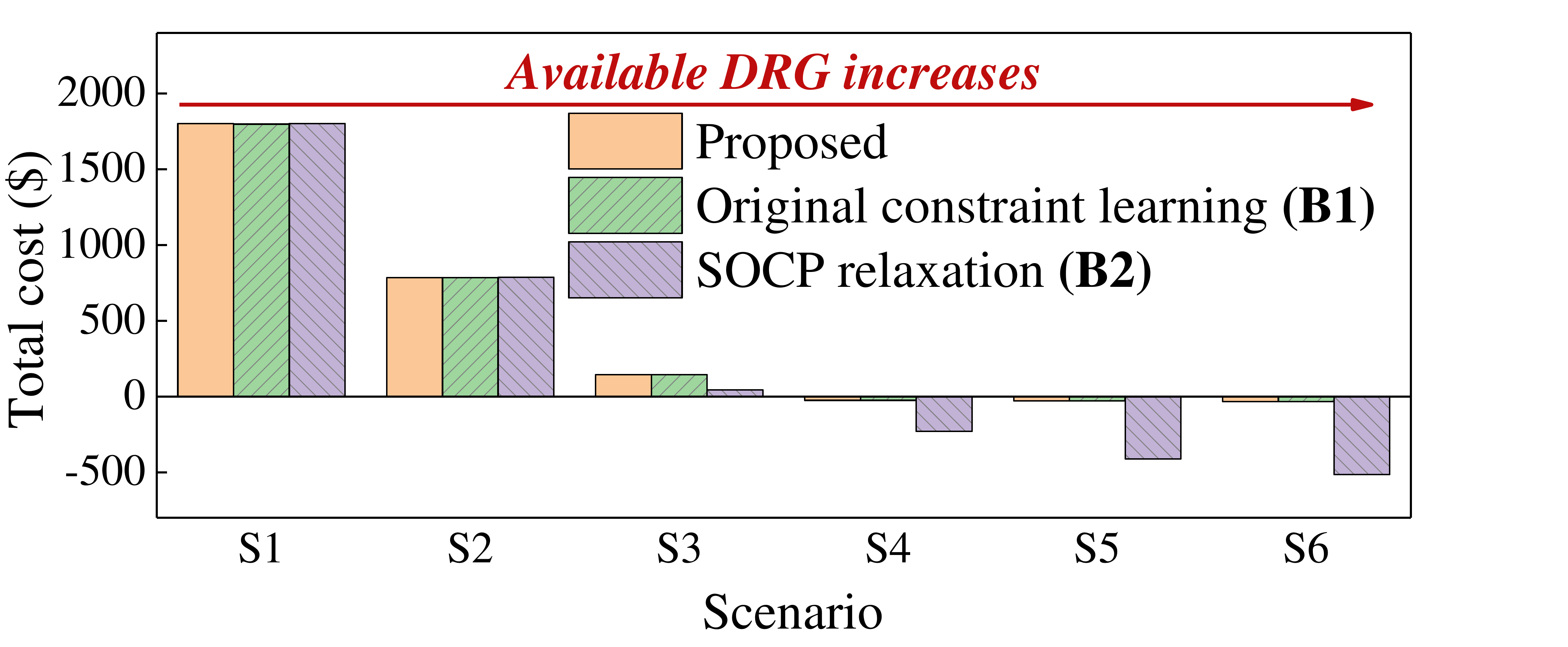}}
	\subfigure[]{\includegraphics[width=0.9\columnwidth]{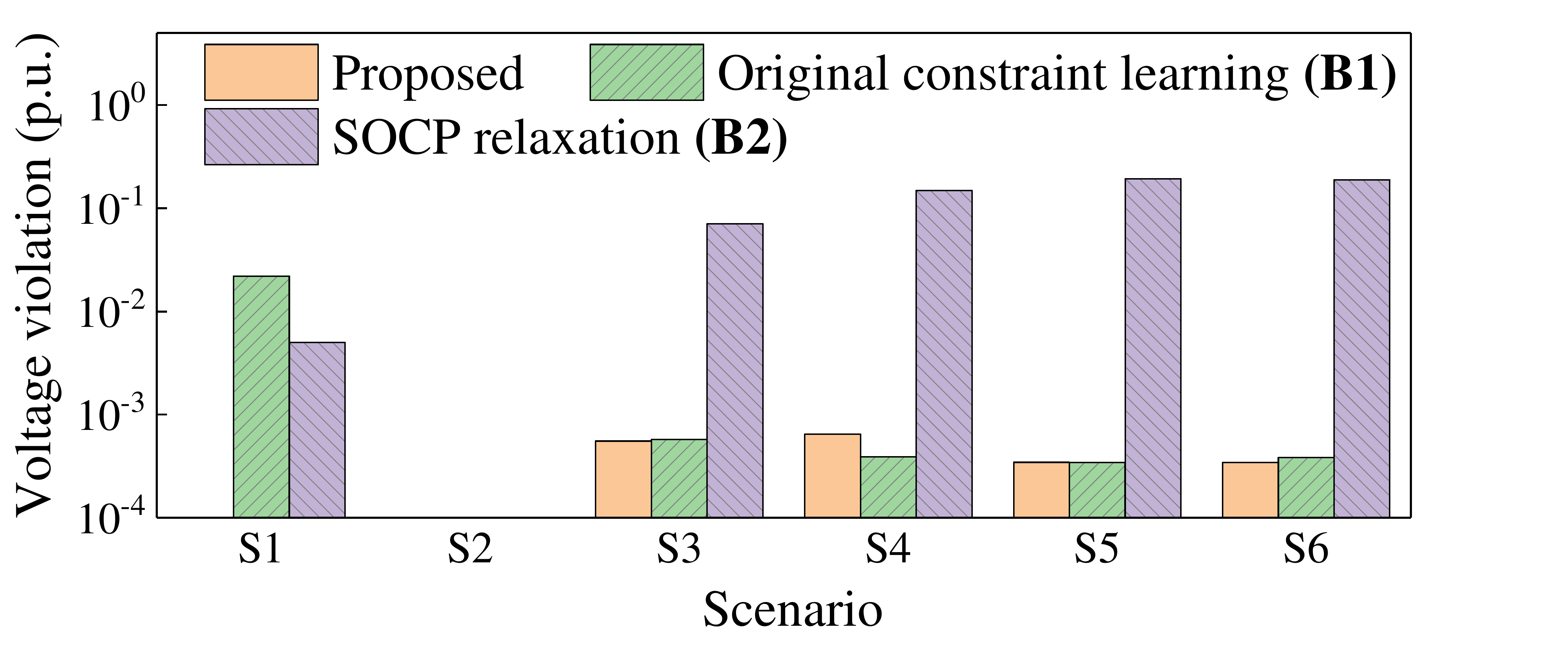}}
	\subfigure[]{\includegraphics[width=0.9\columnwidth]{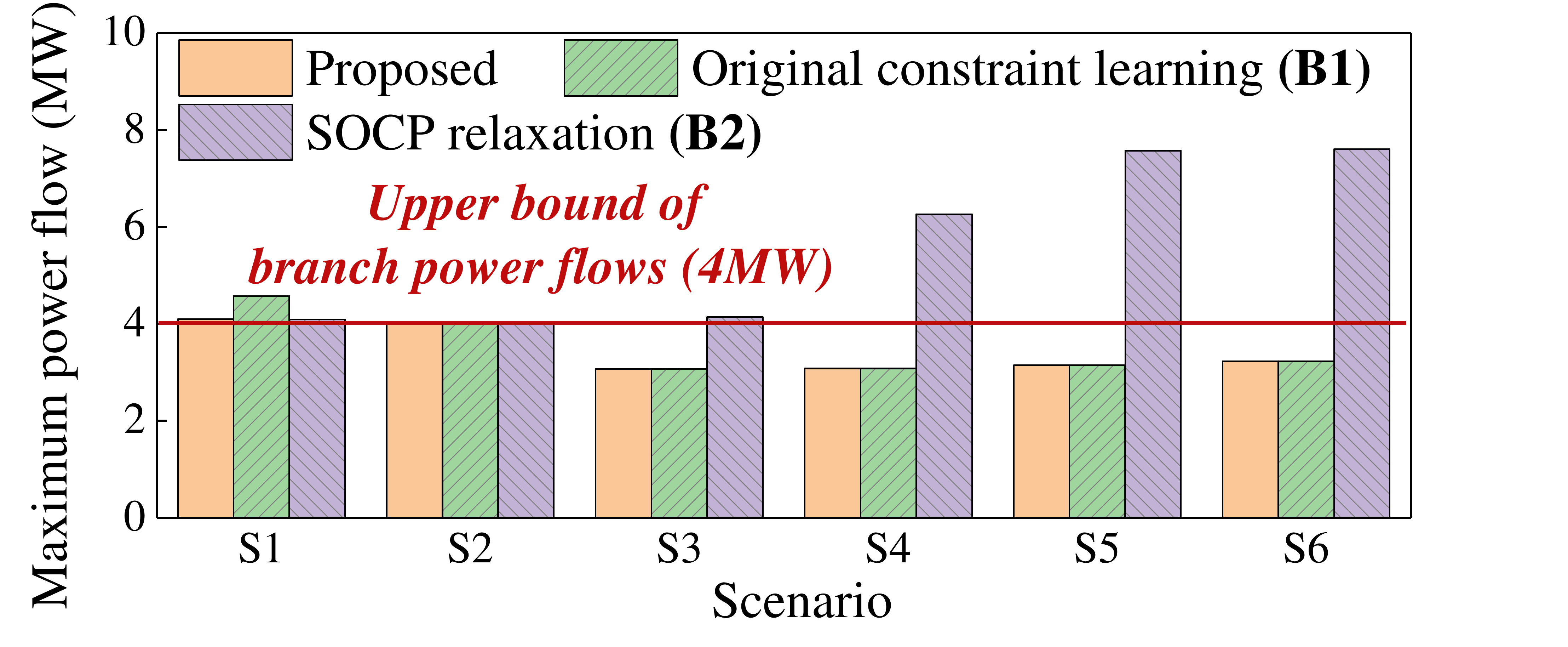}}
 	\caption{Results of (a) total costs, (b) maximum voltage violations, and (c) maximum apparent power flows obtained by different models in the 33-bus case. In (b), the y-axis is expressed in a logarithmic fashion.}
	\label{fig_result_33Bus}
	\vspace{-4mm}
\end{figure}

Fig. \ref{fig_result_33Bus} summarizes the total costs, maximum voltage violations, and maximum apparent power flows of all models in different scenarios. In scenarios S1-S2 with small available DRG, there is no significant reverse power flow. Thus, the results of the SOCP relaxation (\textbf{B2}) have negligible constraint violations and can be regarded as ideal solutions. 
In scenarios S3-S6 with large available DRG, obvious reverse power flows occur in the system, so the SOCP relaxation may be inexact \cite{6815671}. Therefore, its voltage and power flow violations become significant. For example, its maximum voltage violation is around 0.2p.u., while its maximum power flow violation reaches 3.6MW. The conventional constraint learning (\textbf{B1}) can achieve comparable optimality compared to \textbf{B2} in scenarios S1-S2. Moreover, its feasibility is much better than that of \textbf{B2}. For instance, its maximum voltage violation in S2-S6 is less than 0.001p.u.. However, in scenario S1, its maximum violation of bus voltage and power flow limitations reaches 0.02p.u. and 0.57MW, respectively. As mentioned in Section \ref{sec_intro}, these violations are caused by the significant extrapolation from the training data, leading to an unreliable solution. Since the proposed method removes many activation regions from the union in (\ref{eqn_union}), it is an inner approximation of \textbf{B1}. Thus, its total cost is slightly higher than that of \textbf{B1}.
Nevertheless, the maximum relative cost difference between these two methods is always smaller than 0.2\%. Moreover, the proposed method restricts the solution within the activation regions that the Vio-MLP has well learned. Thus, its solution is reliable and shows better feasibility compared to \textbf{B1}. For example, its maximum violations of bus voltage and power flow limitations are only 0.0006p.u. and 0.09MW, which are significantly smaller than those of \textbf{B1}. These results confirm that the proposed method can achieve better feasibility with comparable optimality compared to the conventional constraint learning method.

\subsubsection{Computational efficiency}
\begin{figure}
		\vspace{-4mm}
	\subfigbottomskip=-4pt
	\subfigcapskip=-4pt
	\centering
	{\includegraphics[width=0.9\columnwidth]{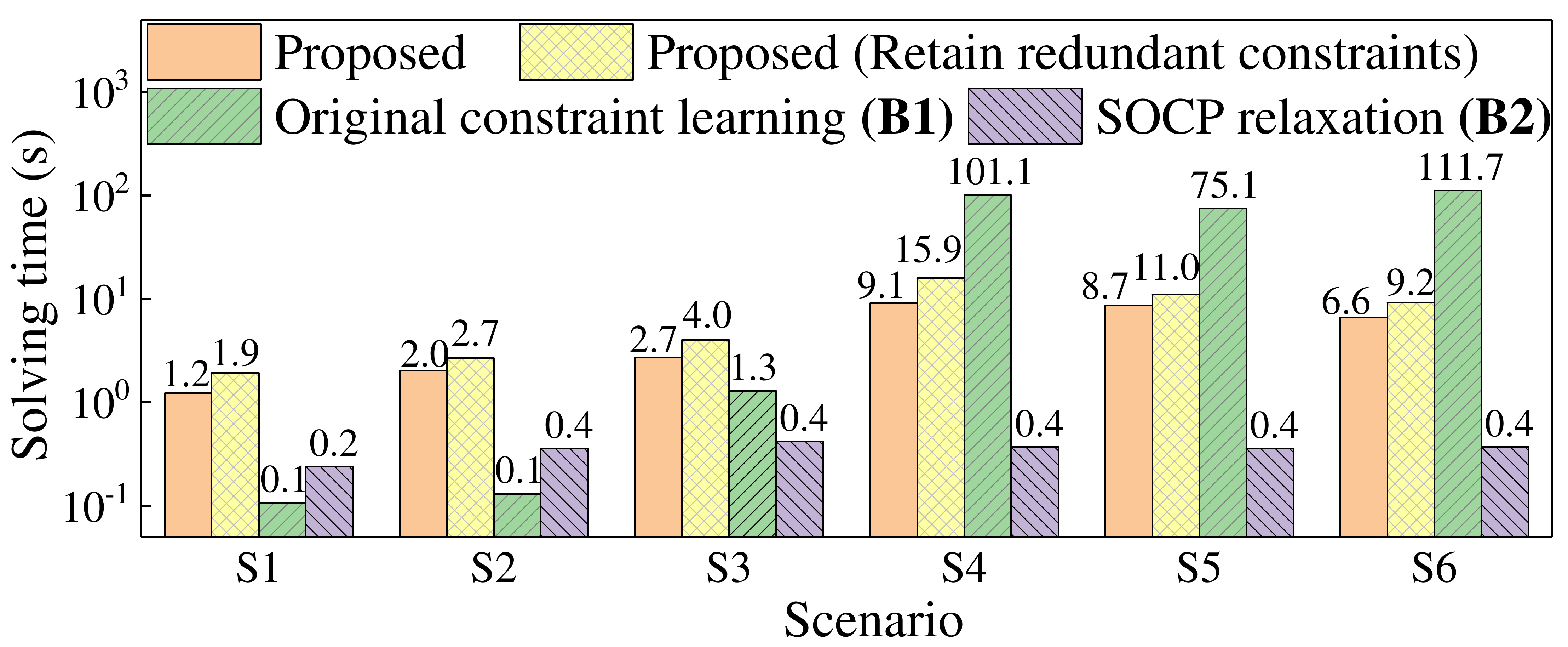}}
	\vspace{-4mm}
 	\caption{Solving times of different models in the 33-bus case. Note the y-axis is expressed in a logarithmic fashion. Here ``Proposed" refers to the proposed model after applying the two-step simplification method. Model ``Proposed (Retain redundant constraints)" represents a new benchmark. This benchmark only drops unnecessary polytopes but retains redundant constraints.}
	\label{fig_solvingTime}
\end{figure}

Fig. \ref{fig_solvingTime} demonstrates the solving time of different models. Here, we introduce a new benchmark called ``Proposed (Retain redundant constraints)" to verify the benefits of removing redundant constraints. In this new benchmark, only the unnecessary polytopes are removed, but the redundant constraints are retained.
The SOCP relaxation \textbf{B2} always achieves great computational efficiency because it is a convex problem. However, it may derive infeasible solutions. 
Compared to the benchmark ``Proposed (Retain redundant constraints)", the proposed method always shows better computational efficiency. This is because it further removes those redundant linear constraints, resulting in a simpler formulation. Although the proposed method spends more time compared to the conventional constraint learning \textbf{B1} in scenarios S1-S3, its time efficiency significantly outperforms \textbf{B1} in the rest scenarios. In scenarios S1-S3, the available DRG is small, so most DRG is consumed locally. Thus, the variation range of the actually used DRG $\bm p_t^\text{DG}$ is small. Meanwhile, the HVAC schedules, i.e., $p_t^\text{HV}$ and $p_t^\text{HV}$, are also restricted in a small range because it is unnecessary to significantly change these schedules to promote the DRG integration. Hence, the feasible region of $\bm x_t$ is tiny and can only lie in a few activation regions. In this case, the majority of activation regions can be directly removed by the pre-solving of the state-of-art solvers in advance. As a result, \textbf{B1} can quickly find the optimal solution. Compared to \textbf{B1}, the proposed method introduces additional linear constraints according to (\ref{eqn_redundant_prune}), so it needs more time to solve each relaxed problem generated in the Branch-and-Bound. Therefore, its solving time is larger than that of \textbf{B1} in S1-S3. In scenarios S4-S6, the variation range of $\bm p_t^\text{DG}$ is large. Meanwhile, considerable DRG is delivered to other nodes, and the corresponding HVAC schedules need to be changed greatly so that more DRG can be integrated. Then, the feasible region of $\bm x_t$ is large and may be distributed in many activation regions. In this case, many activation regions may remain after applying the pre-solving, so it takes much more time for \textbf{B1} to find the optimal solution. Conversely, as mentioned in Section \ref{sec_prune}, the proposed method can remove much more activation regions compared to \textbf{B1}, so its computational efficiency can be guaranteed at a desirable level even in scenarios S4-S6 with large-capacity available DRG. For example, the solving process of the proposed method can be completed in 9.1s in all scenarios, while it may take more than 111s in \textbf{B1}. These results confirm the effectiveness of the proposed simplification method in improving computational efficiency.

\subsection{Case study based on the IEEE 123-bus system} \label{sec_123Bus}

\subsubsection{Parameter setting} 
\begin{figure}
		\vspace{-4mm}
	\centering
	{\includegraphics[width=1\columnwidth]{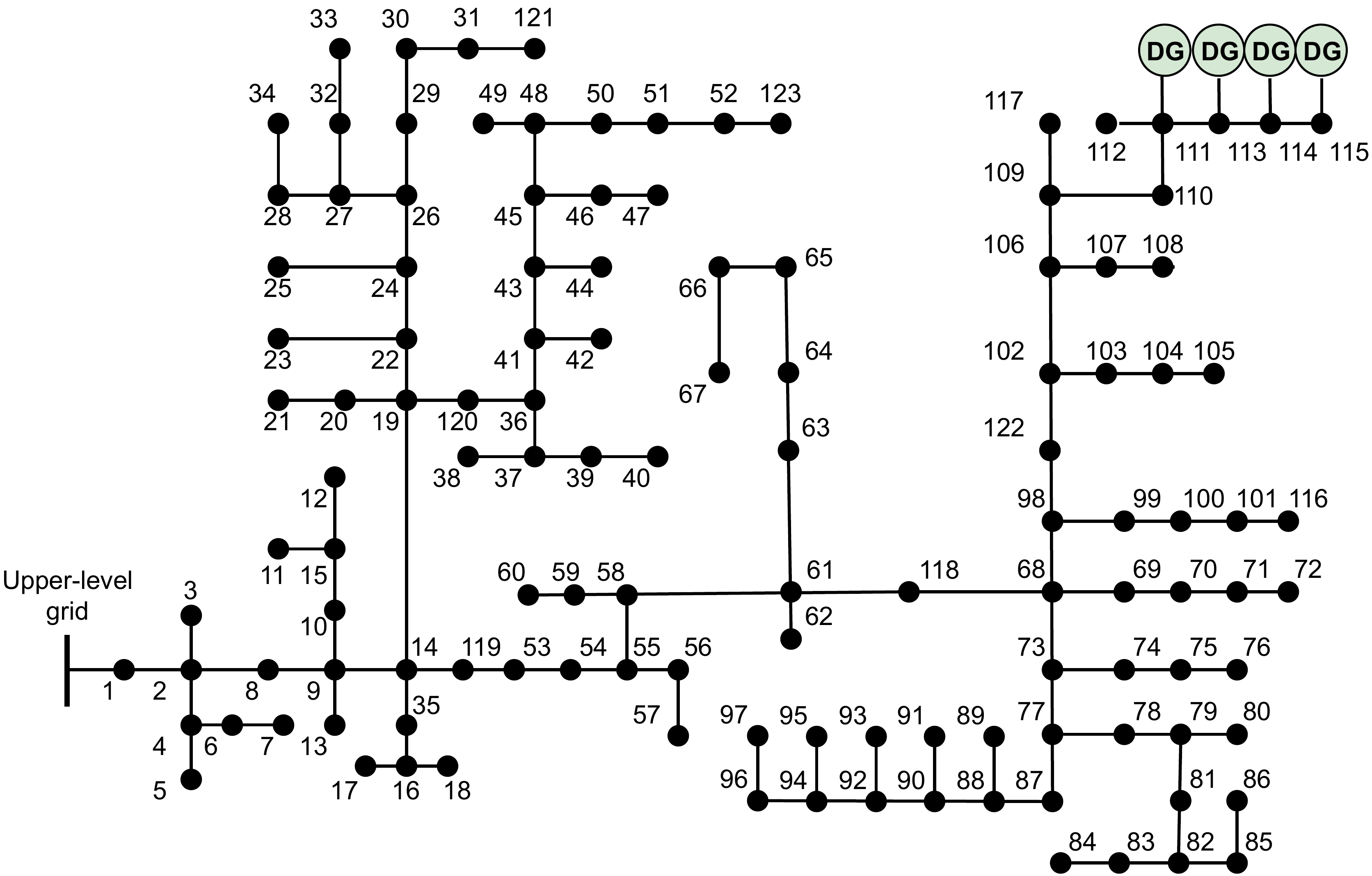}} 
	\vspace{-8mm}
 	\caption{Structure of the IEEE 123-bus system.  }
	\label{fig_123Bus}
	\vspace{-4mm}
\end{figure}

\begin{figure}
		\vspace{-4mm}
	\subfigbottomskip=-4pt
	\subfigcapskip=-4pt
	\centering
	\subfigure[]{\includegraphics[width=0.49\columnwidth]{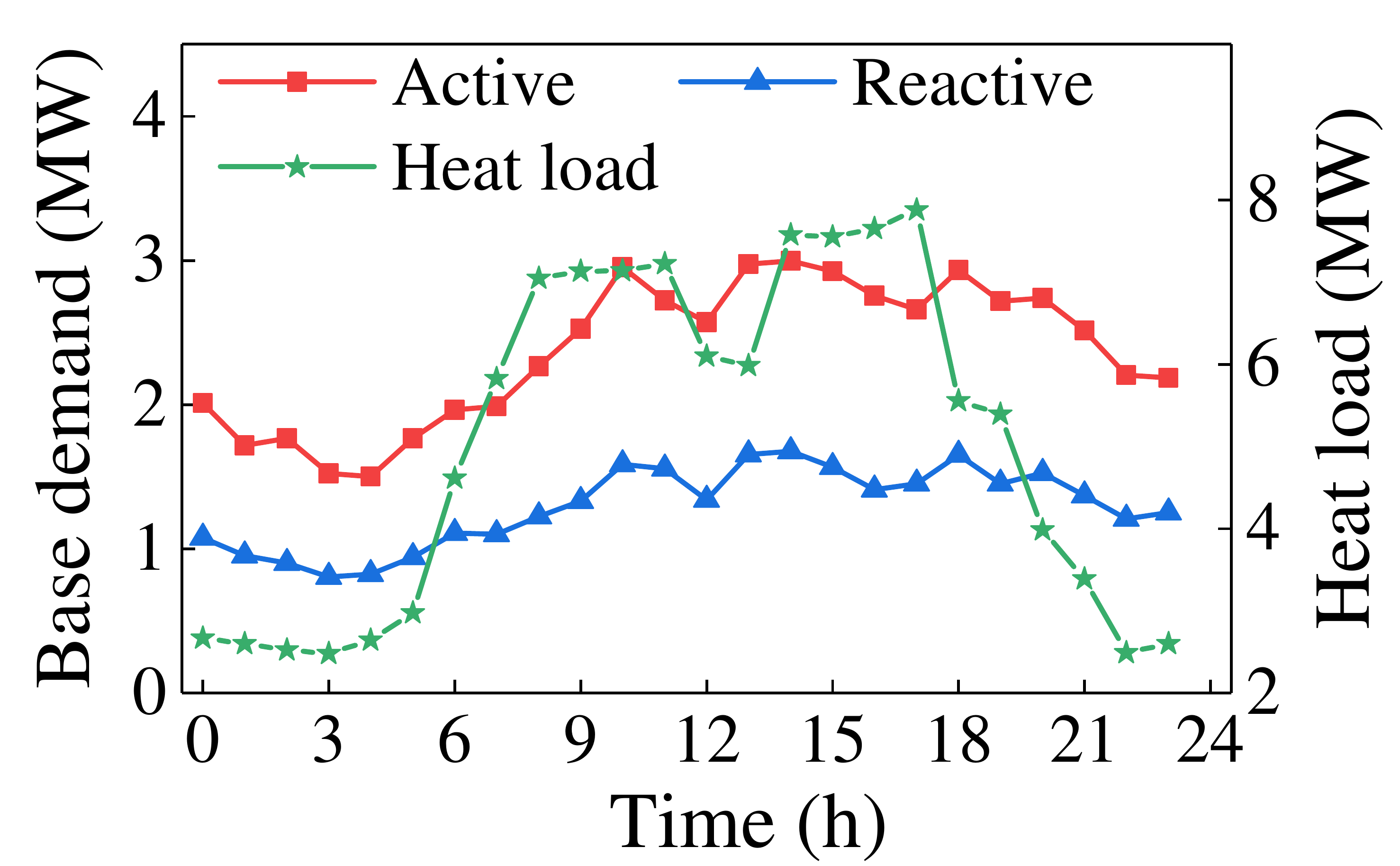}} 
	\subfigure[]{\includegraphics[width=0.49\columnwidth]{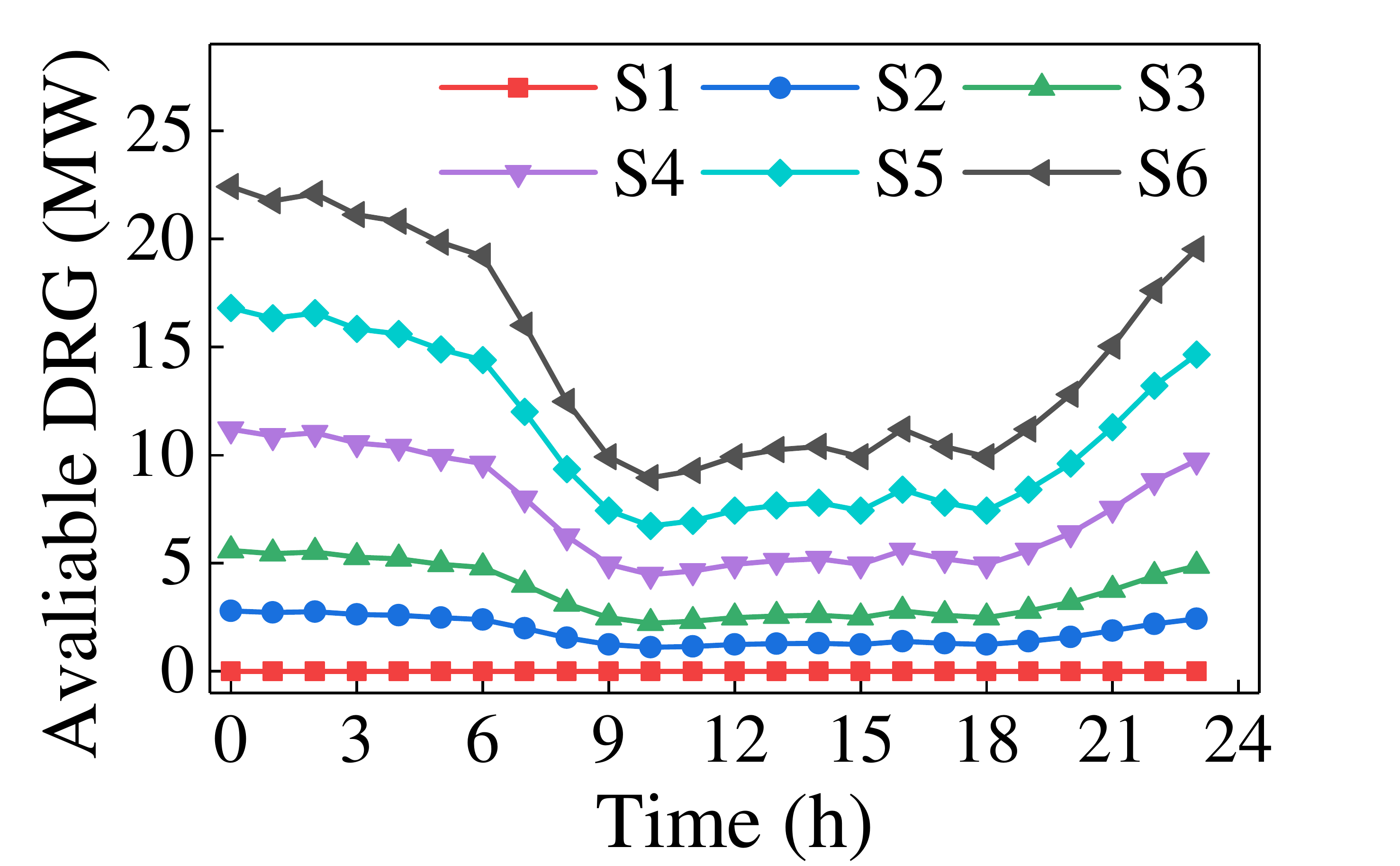}}
 	\caption{Parameters of (a) total active/reactive base power demands and heat load contributed by indoor sources and (b) six scenarios with different available DRG in the 123-Bus case.}
	\label{fig_case_study_123Bus}
	\vspace{-4mm}
\end{figure}

Fig. \ref{fig_123Bus} illustrates the structure of the 123-bus test system equipped with four distributed generators (marked as DG). The voltage at the root node is 4.16kV. The safe region of bus voltages are [0.9p.u., 1.1p.u.], while the maximum allowable branch apparent power flow is 6MW. The unit prices for electricity purchasing/selling and outdoor temperature are the same as those of the 33-bus case. The heat transfer coefficient $g_i$ is 0.02MW/{\textcelsius}. The base electricity demands and indoor heat loads are shown in Fig. \ref{fig_case_study_123Bus}(a). Similarly, we also implement six scenarios with different available DRG, as demonstrated in Fig. \ref{fig_case_study_123Bus}(b).

\subsubsection{Performance of MLPs}
Fig. \ref{fig_MLP_performance} shows the prediction accuracy of the Vio-MLP and Loss-MLP in the 123-bus case, where the structure of the two MLPs are (6, 6, 6) and (4, 4, 4), respectively. Similarly, all sample points approach the red lines, which indicates the excellent performance of MLPs. 
\begin{figure}
		\vspace{-4mm}
	\subfigbottomskip=-6pt
	\subfigcapskip=-6pt
	\centering
	\subfigure[]{\includegraphics[width=0.5\columnwidth]{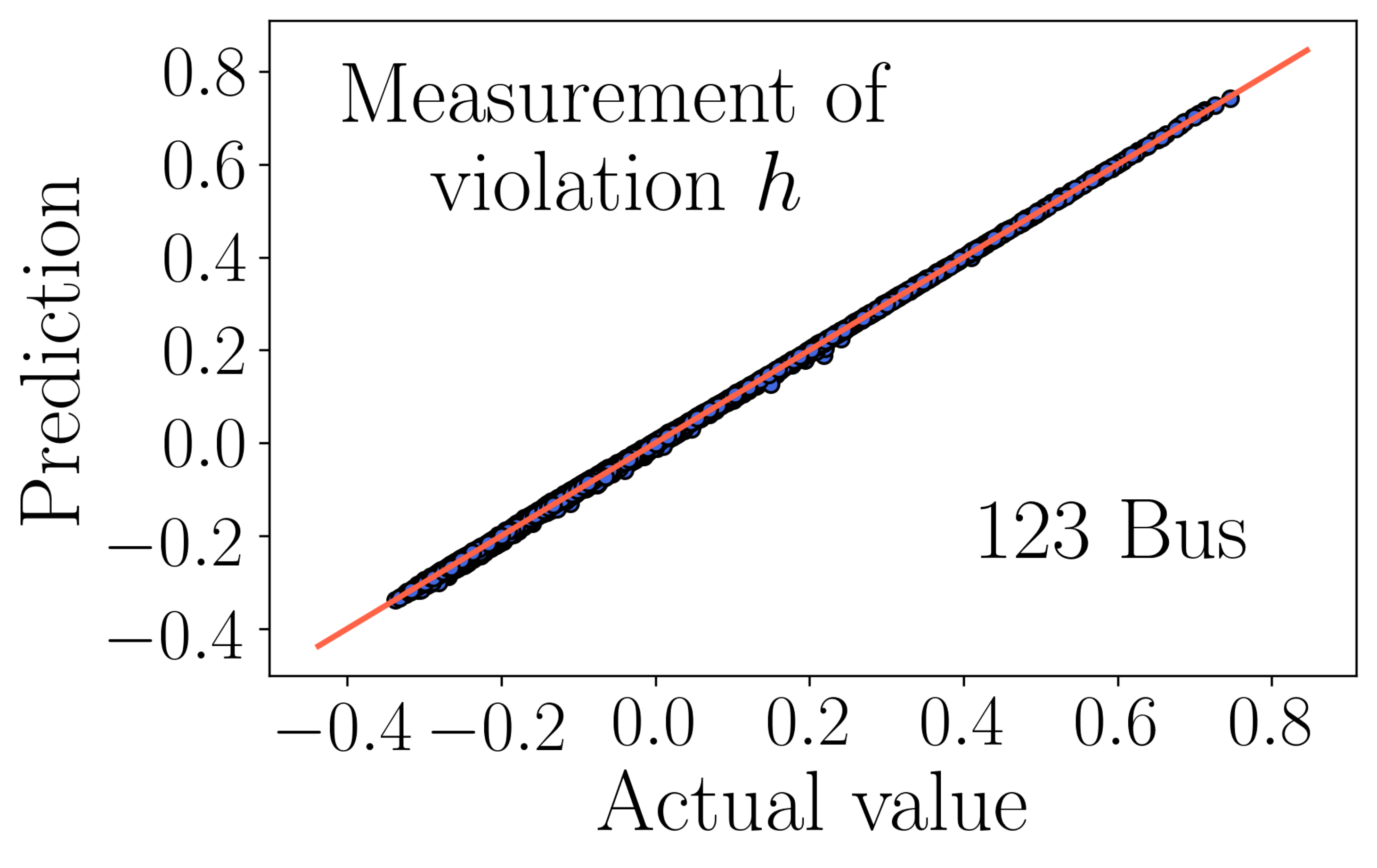}}
	\subfigure[]{\includegraphics[width=0.48\columnwidth]{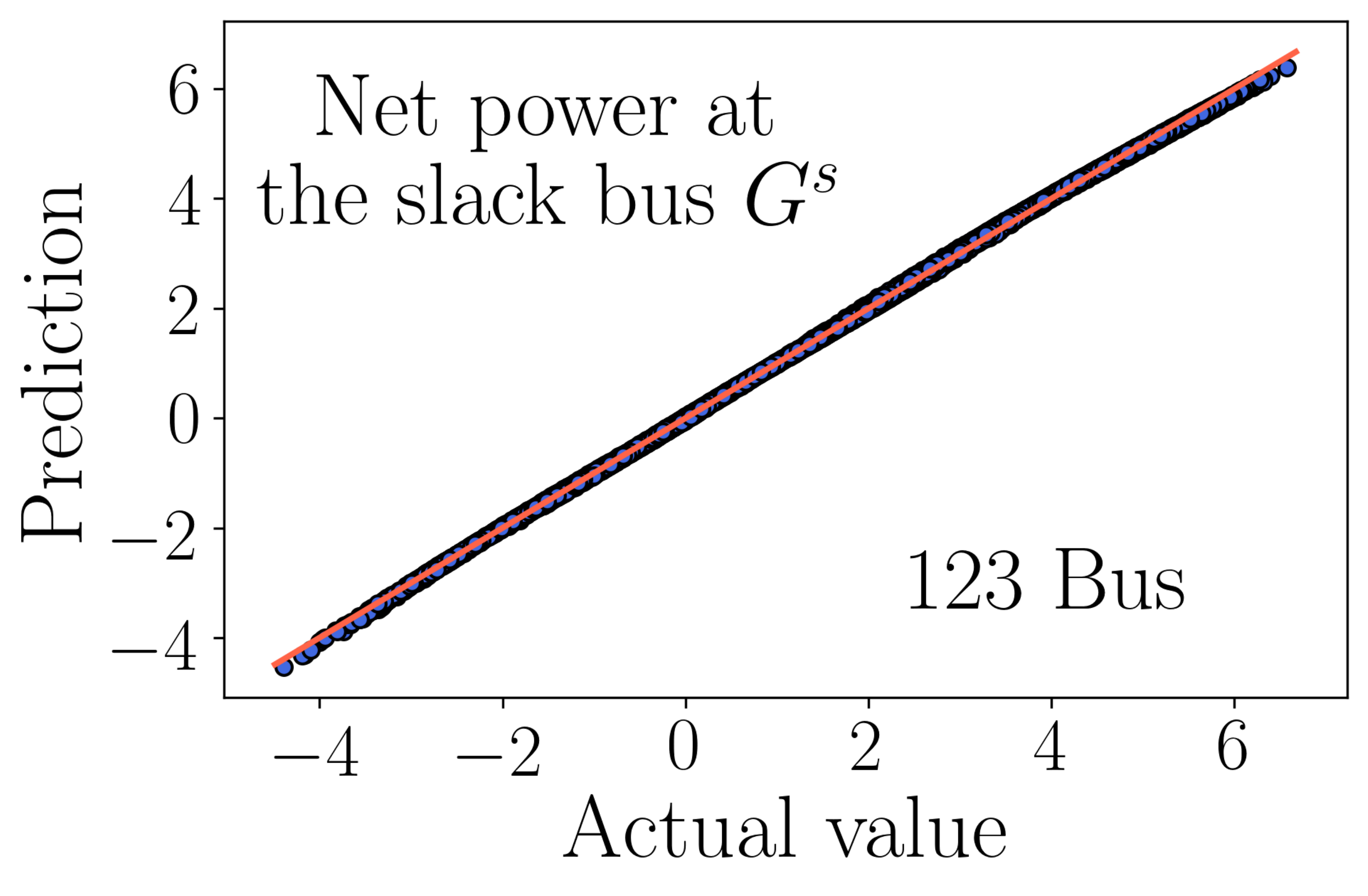}}
	\vspace{-1mm}
 	\caption{Actual and predicted values of (a) the measurement of violations $h_t$ and (b) net power at the substation $G_t^\text{root}$ in the 123-bus case.}
	\label{fig_MLP_performance_123Bus}
	\vspace{0mm}
\end{figure}

\subsubsection{Optimality and feasibility}
\begin{figure}
		\vspace{-4mm}
	\subfigbottomskip=-7pt
	\subfigcapskip=-7pt
	\centering
	\subfigure[]{\includegraphics[width=0.9\columnwidth]{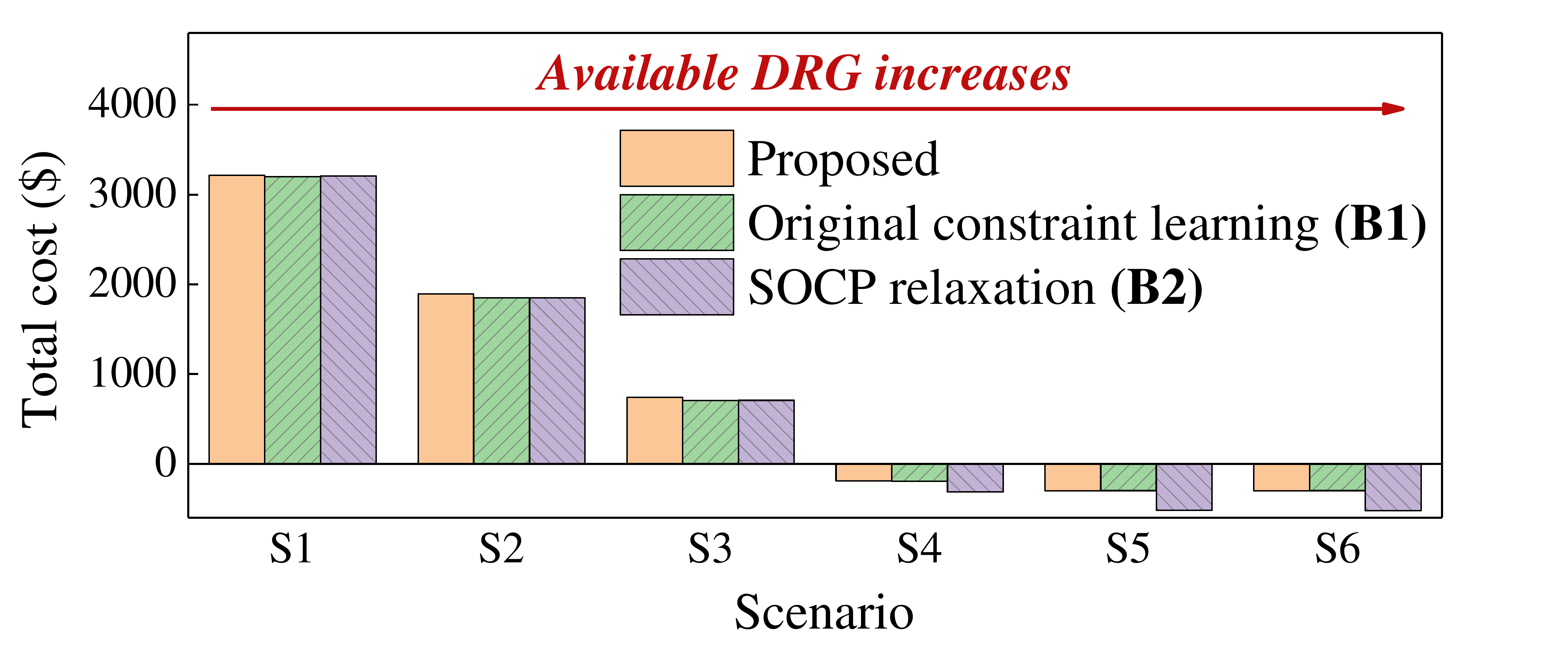}} 
	\subfigure[]{\includegraphics[width=0.9\columnwidth]{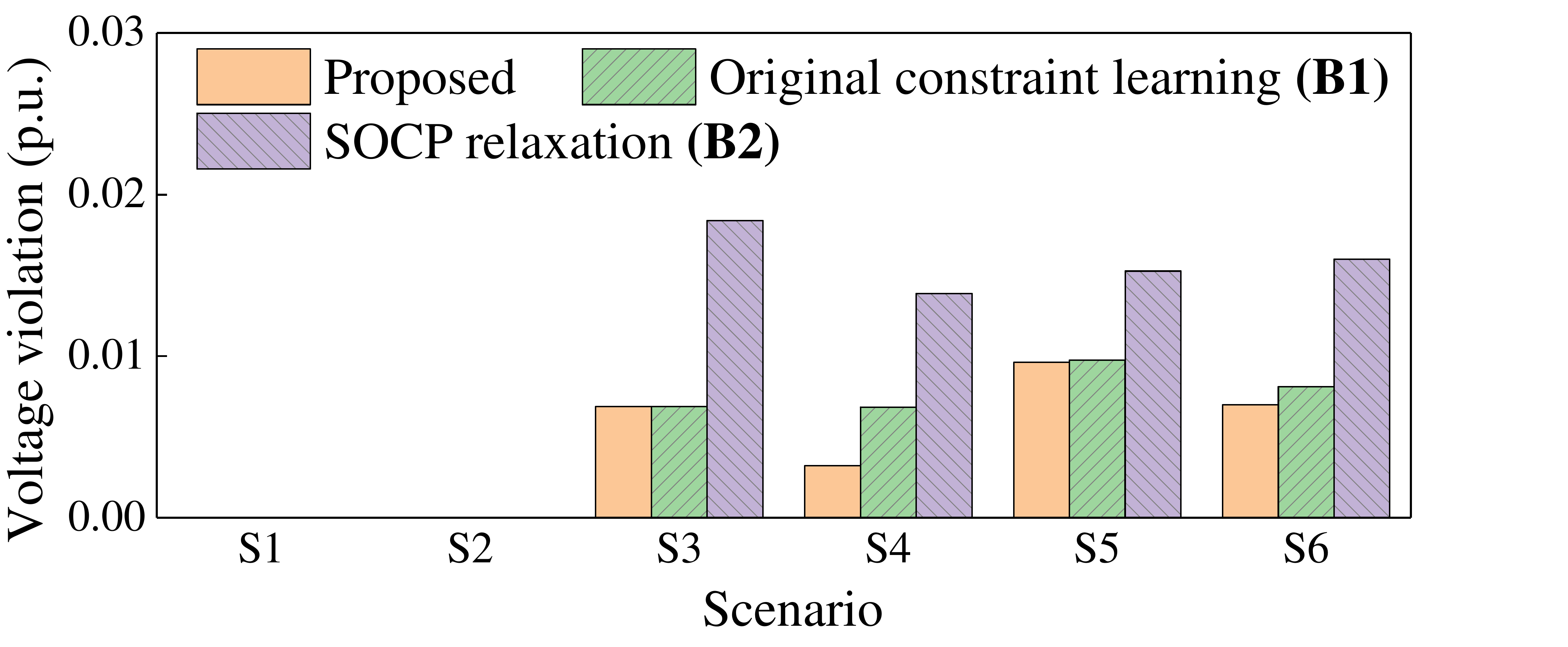}}
	\subfigure[]{\includegraphics[width=0.9\columnwidth]{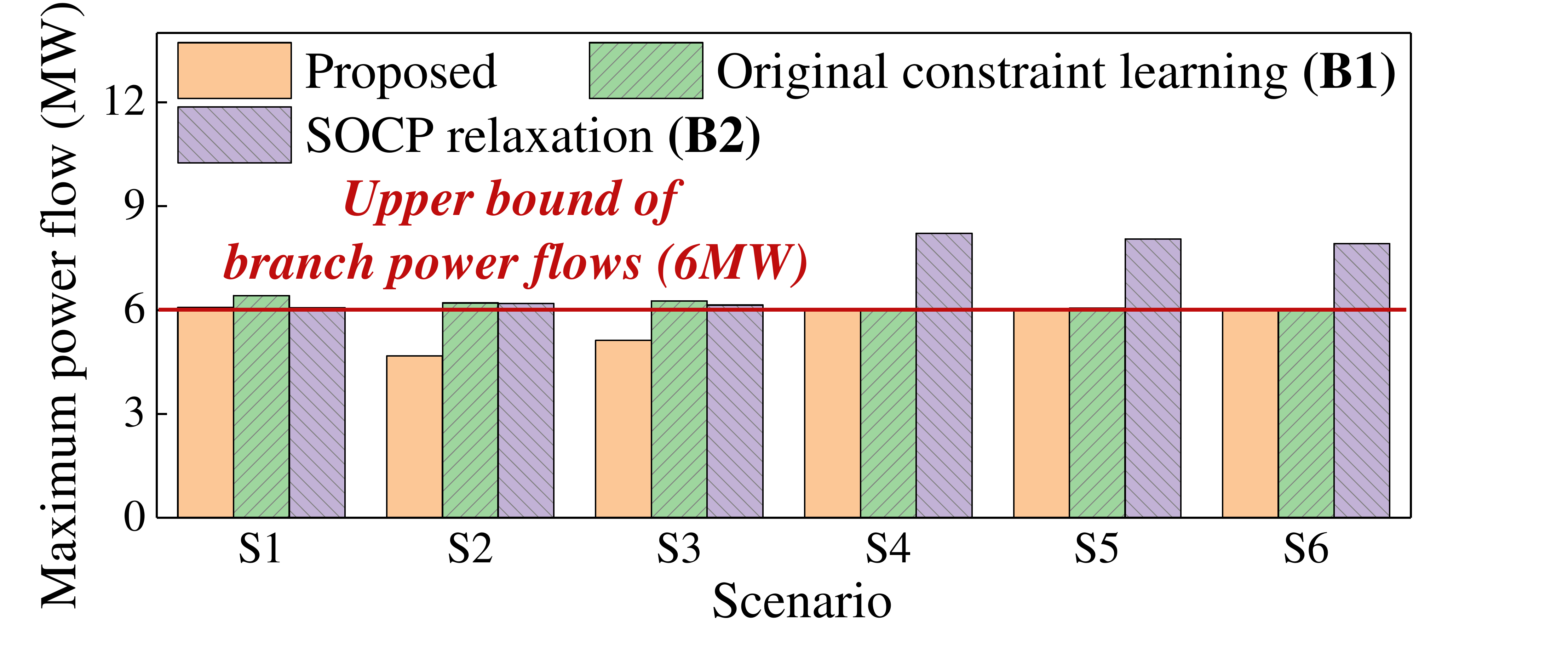}}
 	\caption{Results of (a) total costs, (b) maximum violations of the bus voltage limitation, and (c) maximum apparent power flows obtained by different models in the 123-bus case.}
	\label{fig_result_123Bus}
	\vspace{-4mm}
\end{figure}

Fig. \ref{fig_result_123Bus} summarizes the total costs, maximum voltage violations, and maximum apparent power flows of all models in different scenarios. Similar to the 33-bus case, \textbf{B2} obtains ideal solutions with negligible constraint violations in scenarios S1-S2. However, the violations become obvious in the rest scenarios because the large available DRG causes reverse power flows. The total costs of the conventional constraint learning, i.e., \textbf{B1}, are very close to those of \textbf{B2} in scenarios S1-S2. Moreover, in all scenarios, its voltage violations are smaller than 0.01p.u.. However, its maximum power flow reaches 6.41MW in the first scenario, which is significantly higher than the maximum allowable value (i.e., 6MW). The total cost of the proposed method is only slightly higher than that of \textbf{B1}, but it significantly outperforms \textbf{B1} in feasibility. For example, its maximum power flow violation is only 0.07MW, while it is 0.41MW for \textbf{B1}. These results further validate that the proposed method can achieve near optimality but much better feasibility compared to the conventional constraint learning method.

\subsubsection{Computational efficiency} 
\begin{figure}
	\subfigbottomskip=-4pt
	\subfigcapskip=-4pt
	\centering
	{\includegraphics[width=0.9\columnwidth]{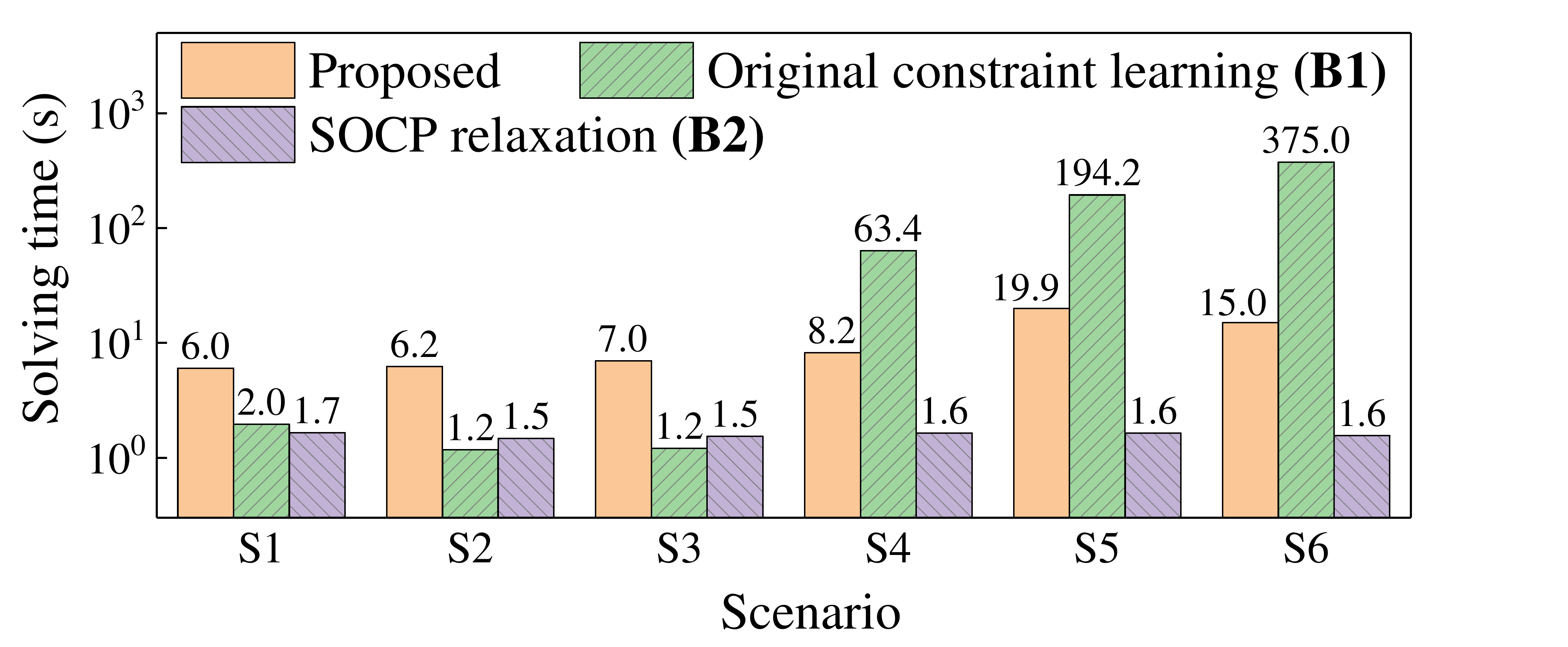}}
	\vspace{-2mm}
 	\caption{Solving times of by different models in the 123-bus case, where the y-axis is expressed in a logarithmic fashion.}
	\label{fig_solvingTime_123Bus}
\end{figure}

Fig. \ref{fig_solvingTime_123Bus} demonstrates the solving time of different models. Similar to Fig. \ref{fig_solvingTime}, the SOCP relaxation \textbf{B2} shows excellent computational efficiency. However, it may lead to infeasible solutions. Compared to the conventional constraint learning \textbf{B1}, the proposed method takes more time in scenarios S1-S3 because it introduces additional linear constraints. Nevertheless, it demonstrates much better computational efficiency in scenarios S4-S6 because the proposed simplification method can remove more activation regions. These results further confirm the desirable computational efficiency of the proposed method.

\section{Conclusions} \label{sec_conclusion}

This paper proposes an efficient constraint learning method to operate ADNs. The proposed method replicates the power flow constraints by training two MLPs with no need for network parameters. A PWL-based interpretation is also presented to explain why the proposed method can accurately replicate complex constraints. Then, a two-step simplification method is proposed to reduce the computational burden of constraint learning. The first step removes unnecessary activation regions, while the second step drops redundant linear constraints. Numerical experiments based on the IEEE 33- and 123-bus systems confirm that the proposed method can achieve better feasibility with comparable optimality compared to the conventional one. Moreover, its computational performance can always be guaranteed at a desirable level.


\footnotesize
\bibliographystyle{ieeetr}
\bibliography{ref}
\end{document}